\documentclass{amsart}

\usepackage{amssymb}
\usepackage{bm}
\usepackage{graphicx}
\usepackage[centertags]{amsmath}
\usepackage{amsfonts}
\usepackage{amsthm}

\newtheorem{Theorem}{Theorem}[section]
\newtheorem{Corollary}{Corollary}[section]
\newtheorem{Lemma}{Lemma}[section]

\newtheorem{Proposition}{Proposition}[section]
\newtheorem{Definition}{Definition}[section]

\newtheorem{Remark}{Remark}[section]

\newtheorem{Notation}{Notation}[section]

\numberwithin{equation}{section}
\newcommand{\beg}{\begin{proof}}
\newcommand{\eproof}{\end{proof}}

 \DeclareMathOperator{\N}{\mathbb{N}}
 
 \DeclareMathOperator{\Q}{\mathbb{Q}}
 \DeclareMathOperator{\R}{\mathbb{R}}
 
 \DeclareMathOperator{\supp}{supp}

 \DeclareMathOperator{\ran}{ran}

\begin{document}

\title{A $c_0$ saturated Banach space with
tight structure} \dedicatory{Dedicated to the memory of Nigel J.
Kalton}
\author{Spiros A. Argyros and Giorgos Petsoulas}
\address{National Technical University of Athens, Faculty of Applied
Sciences, Department of Mathematics, Zografou Campus, 157 80,
Athens, Greece} \email{sargyros@math.ntua.gr, gpetsoulas@yahoo.gr}

\begin{abstract}
It is shown that variants of the HI methods could yield objects
closely connected to the classical Banach spaces. Thus we present
a new $c_0$ saturated space, denoted as $\mathfrak{X}_0$, with
rather tight structure. The space $\mathfrak{X}_0$ is not embedded
into a space with an unconditional basis and its complemented
subspaces have the following structure. Everyone is either of type
I, namely, contains an isomorph of $\mathfrak{X}_0$ itself or else
is isomorphic to a subspace of $c_0$ (type II). Furthermore for
any analytic decomposition of $\mathfrak{X}_0$ into two subspaces
one is of type I and the other is of type II. The operators of
$\mathfrak{X}_0$ share common features with those of HI spaces.
\end{abstract}
\maketitle
\section*{\textbf{Introduction}} The aim of the present paper is
to provide a new norm on $c_{00}(\N)$ resulting a $c_0$ saturated
Banach space. This norm is defined with the use of a modification
of the standard method yielding Hereditarily Indecomposable (HI)
Banach spaces. This approach reveals a Banach space which is $c_0$
saturated but also has a rather tight structure. The following
describes the main properties of the space.
\vspace{2mm}\\
\textbf{Theorem A}
There exists a separable Banach space $\mathfrak{X}_0$
satisfying the following properties.
\begin{enumerate}
\item [(i)] The space $\mathfrak{X}_0$ is $c_0$ saturated and it
is not embedded into a space with an unconditional basis. \item
[(ii)] The dual space $\mathfrak{X}_0^*$ is separable. \item
[(iii)] Every complemented subspace $Y$ of $\mathfrak{X}_0$ is of
one of the following two types. Either $\mathfrak{X}_0$ is
isomorphic to a subspace of $Y$ (type I) or $Y$ is isomorphic to a
subspace of $c_0$ (type II). \item [(iv)] If
$\mathfrak{X}_0=Y\bigoplus Z$ with $Y,Z$ of infinite dimension,
then one of $Y,Z$ is of type I and the other is of type II.
Moreover if the type II complemented subspace is isomorphic to
$c_0$, then the other one is isomorphic to $\mathfrak{X}_0$. In
particular $\mathfrak{X}_0$ is not isomorphic to its square
$\mathfrak{X}_0 \bigoplus \mathfrak{X}_0$.
\end{enumerate}
\vspace{2mm}

Note that properties (iii) and (iv) reminds the strictly
quasi-prime spaces introduced in   \cite{AR}. It is open if the
space $\mathfrak{X}_0$ is strictly quasi-prime. The difference
between $\mathfrak{X}_0$ and the examples of strictly-quasi prime
spaces presented in  \cite{AR} is that the later spaces are not
$c_0$ or $\ell^p$ saturated.

The definition of the norm of $\mathfrak{X}_0$ goes as follows.
We fix two appropriate increasing sequences $(m_j)_{j \in \N}$,
$(n_j)_{j \in \N}$ of natural numbers.
In the first stage we define a norming set $G_0$ as follows.
The set $G_0$ is the minimal subset of $c_{00}(\N)$
satisfying the following properties:
\begin{enumerate}
\item It contains the natural basis $(e_n)_{n \in \N}$ of
$c_{00}(\N)$ and it is symmetric. \item It is closed under the
even operations $(\mathcal{A}_{n_{2j}}, \frac{1}{m_{2j}})$
operation for every $j \in \N$.\\ We recall that this means that
for every $j \in \N,\ d \le n_{2j}$ and $f_1<\ldots<f_d$ in $G_0$
the functional  $\frac{1}{m_{2j}}\sum\limits_{i=1}^{d} f_i$
belongs to $G_0$. Also, as usual, we set $w(f)=m_{2j}$ (the weight
of $f$) if $f$ is a result of a $(\mathcal{A}_{n_{2j}},
\frac{1}{m_{2j}})$ operation. \item $G_0$ contains all
$f=\sum\limits_{i=1}^{n} a_i f_i$, where $\sum\limits_{i=1}^{n}
a_i^2 \le 1$ and $\{f_i\}_{i=1}^{n}$ elements of $G_0$ with
pairwise different weights.
\end{enumerate}

The norm induced by $G_0$ on $c_{00}(\N)$ is denoted as
$\|\cdot\|_{G_0}$. Finally we set $\mathfrak{X}_{G_0}$ the
completion of $(c_{00}(\N),\|\cdot\|_{G_0})$. The space
$\mathfrak{X}_{G_0}$, which is a reflexive one with an
unconditional basis, is a variant of E. Odell and Th. Schlumprecht
space (cf \cite{OS1}) having no $\ell^p$ as a spreading model. We
refer the interested reader to  \cite{AMP} or \cite{AKT} for a
further study of spaces with similar properties.

Next we extend the set $G_0$ to $W_0$ which yields the norm of the
space $\mathfrak{X}_0$ as follows. First we consider a coding
function $\sigma$ similar to the one used in the definitions of HI
and related spaces. Using that coding we define the
$\sigma-n_{2j+1}$ special sequences $(f_i)_{i=1}^{n_{2j+1}}$,
where each $f_i$ belongs to $G_0$. Finally we set
\[W_0=G_0 \cup \{E(\frac{1}{m_{2j+1}} \sum\limits_{i=1}^{n_{2j+1}}
f_i):E \mbox{~interval of~} \N \mbox{~and~} (f_i)_{i=1}^{n_{2j+1}}
\mbox{~is a~} \sigma-n_{2j+1} \mbox{~special sequence~}\}.\] The
space $\mathfrak{X}_0$ is the completion of
$(c_{00}(\N),\|\cdot\|_{W_0})$. Let us point out that the main
difference of $\mathfrak{X}_0$ from a standard HI example (for
example the Gowers-Maurey space (cf \cite{GM})) is that here we
use the odd (i.e. conditional) operations only once at the final
step of the definition of $W_0$. The familiar reader will also
observe that the usual definition of a HI space does not use
condition (3) of the definition of $G_0$. Condition (3) is
critical for proving many of the properties of $\mathfrak{X}_0$ as
well as $\mathfrak{L}(\mathfrak{X}_0)$. The particular use of
this, is in Proposition \ref{L.24}. However it is worth noticing
that the variant of the definition not including condition (3) of
$G_0$ also yields a $c_0$ saturated space and we do not know if
this space satisfies the further properties of the space
$\mathfrak{X}_0$.

As $G_0 \subset W_0$ the identity operator $id:\mathfrak{X}_0
\longrightarrow \mathfrak{X}_{G_0}$ is clearly continuous. The
understanding of the behavior of $id$ on the subspaces of
$\mathfrak{X}_0$ is essential for studying the structure of the
space. In this direction we have the following.
\vspace{2mm}\\
\textbf{Proposition B}
\begin{enumerate}
\item The operator $id:\mathfrak{X}_0 \longrightarrow
\mathfrak{X}_{G_0}$ is strictly singular. \item Let $(x_n)_{n \in
\N}$ be a normalized sequence in $\mathfrak{X}_0$ such that
$\lim\limits_{n} \|x_n\|_{G_0}=0$. Then there exists a subsequence
$(x_n)_{n \in L}$ which is equivalent to the $c_0$ basis. \item
Let $Y$ be a subspace of $\mathfrak{X}_0$ such that $id|_{Y}:Y
\longrightarrow \mathfrak{X}_{G_0}$ is compact. Then $Y$ is
isomorphic to a subspace of $c_0$. \item Let $Y,Z$ be infinite
dimensional subspaces of $\mathfrak{X}_0$ such that
$id|_{Y},id|_{Z}$ are not compact operators. Then $d(S_Y,S_Z)=0$.
\end{enumerate}

Properties (1) and (2) of Proposition B yield that
$\mathfrak{X}_0$ is indeed $c_0$ saturated. The proof of property
(3) requires some beautiful and advanced concepts and results due
to Kalton (cf \cite{K}). Indeed we actually show that any subspace
$Y$ such that $id|_{Y}$ is a compact operator satisfies the $c_0$
tree property (see Definition \ref{L.120}), which according to
Kalton (cf \cite{K}, Thm. 3.2) yields that $Y$ is isomorphic to a
subspace of $c_0$. A consequence of properties (3) and (4) is that
every complemented subspace of $\mathfrak{X}_0$ is either of type
I or of type II (see Thm. A).

A second result describing the tight structure of $\mathfrak{X}_0$
and its relation to HI spaces concerns the operators. Let point
out that Sobczyk's theorem (cf \cite{SO}) yields that
$\mathfrak{X}_0$ admits many projections as every $c_0$ subspace
is a complemented one. The following  explains that the non
strictly singular operators on $\mathfrak{X}_0$ have a precise
structure.
\vspace{2mm}\\
\textbf{Theorem C} Every bounded linear operator $T:\mathfrak{X}_0
\longrightarrow \mathfrak{X}_{0}$ is of the form $T=\lambda \cdot
I+S$ with $S$ satisfying the following. If there exists a subspace
$Y$ of $\mathfrak{X}_0$ with $S|_{Y}$ is an isomorphism, then $Y$
is isomorphic to a subspace of $c_0$.

The paper is organized into six sections. In the first one we
present the definition of the norming sets $G_0$ and $W_0$ and the
corresponding spaces $\mathfrak{X}_{G_0}$ and $\mathfrak{X}_0$.
The second one is devoted to the study of bounded block sequences
$(x_n)_{n \in \N}$ in $\mathfrak{X}_0$ satisfying $\overline{\lim}
\ \|x_n\|_{G_0}>0$. The main result of this section (Prop.
\ref{L.25}) asserts that any such sequence contains arbitrarily
large seminormalized averages on which the $G_0$ and $W_0$ norms
coincide. Sections 3 and 4 are devoted to the basic inequality,
the exact pairs and the dependent sequences. All these are closely
related to the corresponding concepts and results appeared and
used in the study of HI spaces (cf
 \cite{AMP}). Sections 5 and 6 include the proofs of Theorems A and C
respectively.

\vspace{2mm} We make use of the following standard notation
throughout
this article.\\
\begin{enumerate}
\item[i.] We denote by $c_{00}(\N)$ the vector space
$c_{00}(\N)=\{f:\N\to\R : f(n)\neq 0 \text{ for finitely many }
n\in\N\}$ and by $c^{\Q}_{00}(\N)$ the set of all elements of
$c_{00}(\N)$ with rational coordinates. For every $x\in
c_{00}(\N)$ we denote by $\supp x$ the set $suppx=\{n\in\N:
x(n)\neq 0\}$ and by $\ran x$ the minimal interval of $\N$ that
contains $\supp x$. \item[ii.] We denote by $(e_n)_n$ the standard
Hamel basis of $c_{00}(\N)$, which will also be considered as
functionals on $c_{00}(\N)$ acting through the usual inner product
and denoted as $(e_n^*)_n$. \item[iii.] Let $E_1,E_2$ be two
nonempty finite subsets of $\N$. We write $E_1<E_2$ if $\max E_1<
\min E_2$. Also for a $n \in \N$, we write $n<E_1$ if $\{n\}<E_1$.
If $x_1,x_2$ are non zero sequences of $c_{00}(\N)$ we write
$x_1<x_2$ whenever $\ran x_1< \ran x_2$. In addition for a
sequence $f:\N\to\R$ and $E$ an interval of $\N$ we denote by $Ef$
the sequence $f\cdot X_E$, where $X_E$ is the characteristic
function of $E$. \item[iv.] We say that a subset $F$ of
$c_{00}(\N)$ is closed under the $(\mathcal{A}_n,
\theta)$-operation for $n \in \N$ and $0<\theta<1,$ if for every
$d \le n$ and for every $f_1<\ldots<f_d$ in $F$ we have that
$\theta\sum_{i=1}^d f_i \in F$. \item [v.] We say that a subset
$F$ of $c_{00}(\N)$ is symmetric, if for every $f \in F$ it
follows that $-f \in F$. \item [vi.] We say that a subset $F$ of
$c_{00}(\N)$ is closed in restrictions to finite intervals of
$\N$, if for every $f \in F$ and $E$ finite interval of $\N$, it
follows that $Ef \in F$. \item[vii.] Let $L$ an infinite subset of
$\N$ and $k \in \N$. We denote by $[L]^k$ the set of all subsets
of $\N$ with $k$ elements and with $[L]$ the set of all infinite
subsets of $L$.
\end{enumerate}

\section{\textbf{The norming set of the Banach space}
$\mathfrak{X}_0$} In this section we define the norming sets $G_0$
and $W_0$ yielding the spaces $\mathfrak{X}_{G_0}$ and
$\mathfrak{X}_0$ respectively. We fix two sequences of natural
numbers $(m_j)_j$ and $(n_j)_j$ defined recursively as follows. We
set $m_1=2^8$ and $m_{j+1}=m_j^5$ and $n_1=2^7$ and
$n_{j+1}=(2n_j)^{s_{j+1}}$ where $s_{j+1}=\log_2(m^4_{j+1})$,
$j\geq 1$.

\begin{Definition}
Let $G_0$ be the minimal subset of $c_{00}(\N)$ satisfying the
following:
\begin{enumerate}
\item $G_0$ contains the set $F_0=\{e^*_n:n\in \N\}$. \item $G_0$
is symmetric. \item $G_0$ is closed under the
$(\mathcal{A}_{n_{2j}}, \frac{1}{m_{2j}})$ operation for every $j
\in \N$. \item It contains the set $\{\sum_{i=1}^d a_i f_i:\;d\in
\N,\ a_i\in\Q,\sum_{i=1}^d a^2_i\leq 1$ and $f_i \in G_0$ with
$(w(f_i))_{i=1}^d$ pairwise different\}.
\end{enumerate}
\end{Definition}
\vspace{3mm} For an $f \in G_0$ we say that $f$ has weight
$m_{2j}$ and we write $w(f)=m_{2j}$ if and only if there exists $d
\in \N$ with $d\le n_{2j}$  and $f_1<\ldots<f_d$ in $G_0$ such
that $f=\frac{1}{m_{2j}}\sum_{i=1}^d f_i$. Such an $f$ is called a
functional with weight.

\begin{Definition}
We define the Banach space $\mathfrak{X}_{G_0}$ to be the
completion of $(c_{00}(\N),\|\cdot\|_{G_0})$, where
\[ \|x\|_{G_0}=\sup\{|f(x)|:\;f\in G_0 \},x\in
c_{00}(\N) .\]
\end{Definition}
\vspace{3mm} The space $\mathfrak{X}_{G_0}$ resembles the space
defined by E. Odell and Th. Schlumprecht in \cite{OS1}, where it
is proved that this space does not admit $c_0$ and $\ell^p$ as a
spreading model. Since the space $\mathfrak{X}_{G_0}$ does not
admit $c_0$ and $\ell^p$ as a spreading model for every $1 \le
p<\infty$ (for a proof we refer to \cite{OS1} or \cite{AMP}) and
the basis is unconditional, it follows that this space is
reflexive. Next we define the  $\sigma-n_{2j+1}$ special sequences
and the norming set $W_0$.

\begin{Definition}\label{L.51}
Let $\Omega_1,\Omega_2$ be two disjoint infinite subsets of $\N$
and $Q_s=\{(f_1,\ldots,f_d):\;d\in \N,f_i \in G_0,f_i \neq
0,i=1,\ldots,d,f_1< \ldots <f_d\}$. Since $Q_s$ is a subset of
$c_{00}^{\Q}(\N)$, it follows that $Q_s$ is countable, so we may
select an injective map $\sigma:Q_s \to \{2j:\;j \in \Omega_2\}$
such that
\[m_{\sigma(f_1,\ldots,f_d)}>
\max\{\frac{1}{|f_i(e_l)|}:\;i=1,...,d,\ l \in \supp(f_i)\} \cdot
\max\supp(f_d)\]
for every $(f_1,\ldots,f_d)$ in $Q_s$. \\
Let $j \in \N$. A finite sequence $(f_i)_{i=1}^{n_{2j+1}}$ with
$(f_1,\ldots,f_{n_{2j+1}})$ in $Q_s,$ is said to be
$\sigma-n_{2j+1}$ special sequence provided:
\begin{enumerate}
\item each $f_i$ is a functional with weight.
\item $w(f_1)=m_{2j_1}$, $j_1\in\Omega_1$ and
$n_{2j+1}^2<m_{2j_1}$.
\item $w(f_{i+1})=m_{\sigma(f_1,...,f_i)}$, for all
$i\in \{1,\ldots,n_{2j+1}-1\}$.
\end{enumerate}
\end{Definition}
\vspace{3mm}
We pass now to define the norming set $W_0$ and the
corresponding space $\mathfrak{X}_0$.

\begin{Definition}
Let $W_0$ be the minimal subset of  $c_{00}(\N)$ such
that
\begin{enumerate}
\item $G_0\cup \{\frac{1}{m_{2j+1}}\sum\limits_{i=1}^{n_{2j+1}}
f_i:\; f_i \in G_0$ with $(f_i)_{i=1}^{n_{2j+1}} \sigma-n_{2j+1}$
special sequence\} $\subset W_0$
\item $W_0$ is symmetric
\item $W_0$ is closed in restrictions to the finite intervals
of $\N$.
\end{enumerate}
The Banach space $\mathfrak{X}_0$ is the completion of
$(c_{00}(\N),\|\cdot\|_{W_0})$, where
\[ \|x\|_{W_0}=\sup\{|f(x)|:\;f\in W_0 \},x\in
c_{00}(\N).\]
\end{Definition}
\vspace{3mm} For an $f \in W_0$ we say that $f$ has weight
$m_{2j+1}$ and we write $w(f)=m_{2j+1}$ if and only if there
exists a $\sigma-n_{2j+1}$ special sequence
$(f_i)_{i=1}^{n_{2j+1}}$ in $G_0$ such that $f=\epsilon
E\frac{1}{m_{2j+1}}\sum_{i=1}^d f_i$, where $|\epsilon|=1$ and $E$
a finite interval of $\N$.

\begin{Remark}
\begin{enumerate}
\item The norming set $G_0$ is closed in restrictions to the
finite subsets of $\N$. \item It is easily checked that
$W_0=G_0\cup \{\epsilon
E\frac{1}{m_{2j+1}}\sum\limits_{i=1}^{n_{2j+1}}
f_i:\;|\epsilon|=1, E$ finite interval of $\N$ and $f_i \in G_0$
with $(f_i)_{i=1}^{n_{2j+1}} \sigma-n_{2j+1}$ special sequence\}.
\item If $f\in W_0,$ then $\|f\|_{\infty}\le 1$. \item The basis
$(e_n)_{n\in \N}$ of the Banach space $\mathfrak{X}_0$ is
bimonotone and $\|e_n\|_{W_0}=1$ for all $n\in\N$. Also the basis
of $\mathfrak{X}_{G_0}$ is 1-unconditional (i.e. for every $x \in
\mathfrak{X}_{G_0}$ and $E$ subset of $\N$ we have that
$\|Ex\|_{G_0}\le \|x\|_{G_0}$).
\end{enumerate}
\end{Remark}

\section{\textbf{Estimating averages in $\mathfrak{X}_0$}}
The main result of this section is the following proposition,
which is a key ingredient for studying the structure of the
space $\mathfrak{X}_0$ and the corresponding one of
$\mathfrak{L}(\mathfrak{X}_0).$
\vspace{2mm}\\
\textbf{Proposition 4.1.}
For every $\epsilon>0$ there exists $n \in \N$
such that for every $k \in \N$ with $k>n$ and every
block sequence $(x_n)_{n \in \N}$ with
\[0<\epsilon<\|x_n\|_{G_0}\le\|x_n\|_{W_0}\le
1 \mbox{~for all~} n \in \N\]
there exists an $L\in [\N]$ such that
\\ for every $n_1<\ldots<n_k$ in $L$ and every
$\phi\in (W_0\backslash G_0)$ it follows that
\[\phi(\frac{x_{n_1}+\ldots+x_{n_k}}{k})\le
\|\frac{x_{n_1}+\ldots+x_{n_k}}{k}\|_{G_0}.\]
Hence
\[\|\frac{x_{n_1}+\ldots+x_{n_k}}{k}\|_{W_0}=
\|\frac{x_{n_1}+\ldots+x_{n_k}}{k}\|_{G_0}.\]
\vspace{2mm}\\
The proof of the proposition is of combinatorial nature, i.e. is
mainly based on Ramsey's theorem \cite{R}, and uses the property
that the norming set $G_0$ is closed under rational convex
combinations. Actually this is the only point where that property
is used. A consequence of the proposition is that every block
sequence $(x_n)_n$ with $\overline{\lim} \|x_n\|_{W_0}< \infty$
and $\underline{\lim} \|x_n\|_{G_0}>0$ admits further block
sequences which are arbitrarily large seminormalized $l^1$
averages and moreover their norm in $\mathfrak{X}_0$ coincide with
the corresponding one in $\mathfrak{X}_{G_0}$. This is the
fundamental ingredient for proving the properties of the space
$\mathfrak{X}_0$ and the properties of
$\mathfrak{L}(\mathfrak{X}_0)$, as it permits to pass to exact
pairs in $\mathfrak{X}_0$ and then to dependent sequences in
$\mathfrak{X}_0$.

\begin{Definition}
Let $k \in \N$, $\delta>0$ and $(x_n)_{n=1}^{k}$ be a finite block
sequence. Let also $\phi$ be a functional in $(W_0 \backslash
G_0)$ of the form $\phi=\frac{1}{m_{2j+1}}\sum\limits_{i=1}^{q}
f_i$, where $q \le n_{2j+1}$ and $f_1<\ldots<f_q$ in $G_0$.
\begin{enumerate}
\item We will say that $(x_n)_{n=1}^{k}$ is $(\phi,\delta)$
separated if there exist
\begin{enumerate}
\item [i.] $d_1<d_2<d_3$ in $\{1,\ldots,k\}$ and
\item [ii.] $E_{d_1}<E_{d_2}<E_{d_3}$ subintervals of
$\{1,\ldots,q\}$
\end{enumerate}
such that $\delta<\sum\limits_{p \in E_{d_i}} f_p(x_{d_i})$
for all $i=1,2,3$.
\item We will say that $(x_n)_{n=1}^{k}$ is $\delta-$
separated if there exists a functional $\phi$ in
$(W_0 \backslash G_0)$ such that
$(x_n)_{n=1}^{k}$ is $(\phi,\delta)$ separated.
\end{enumerate}
\end{Definition}

\begin{Notation}
Let $k \in \N,\ n_1<\ldots<n_k$ in $\N$ and $(x_{n_i})_{i=1}^k$ be
a finite block sequence. Let also $\phi \in (W_0 \backslash G_0)$
with $\phi=\frac{1}{m_{2j+1}} \sum\limits_{i=1}^{q} f_i$, where
$q\le n_{2j+1}$ and $(f_i)_{i=1}^{q}$ successive elements in
$G_0$.
\\ Then for every $d \in \{1,\ldots,k\}$ we define the set
$E_{n_d}^{\phi}$ as follows:
\[ E_{n_d}^{\phi}=\{i \in \{1,\ldots,q\}:ran(f_i)\cap
ran(x_{n_d})\neq \emptyset\}.\]
\end{Notation}

\begin{Lemma}\label{L.64}
Let $k \in \N,\ \delta>0$ and $(x_n)_{n \in \N}$ be a block
sequence. Let also $L \in [\N]$ satisfying that for every
$n_1<\ldots<n_k$ in $L$ there exists $\phi \in (W_0 \backslash
G_0)$ such that the block sequence $(x_{n_i})_{i=1}^{k}$ is
$(\phi,\delta)$ separated. Then for every $M>0$, there exists $L_0
\in [L]$ such that $M\le \|x_n\|_{G_0}$ for all $n \in L_0$.
\end{Lemma}
\noindent{\bf{Proof:}} Assume on the contrary. Then there exist
$M>0$ and $P \in [L]$ (assume without loss of generality that
$P=L$)
such that $\|x_n\|_{G_0}\le M$ for all $n \in L$.\\
From the fact that for every $n_1<\ldots<n_k$ in $L$ there exists
$\phi \in (W_0 \backslash G_0)$ such that the block sequence
$(x_{n_i})_{i=1}^{k}$ is $(\phi,\delta)$ separated, applying
Ramsey's Theorem \cite{R} for the set $[L]^k$, (i.e. if
$[L]^k=A\cup B$, then there exists $L_1 \in [L]$ such that
$[L_1]^k \subset A$ or $[L_1]^k \subset B$), we conclude that
\begin{enumerate}
\item there exists $L_1 \in [L]$ and
\item there exist $d_1<d_2<d_3$ in $\{1,\ldots,k\}$
\end{enumerate}
such that for every $F=\{n_1<\ldots<n_k\} \in [L_1]^k$,
there exists $\phi_F \in (W_0\backslash G_0)$ of the form
$\phi_F=\frac{1}{m_{2j_F+1}}\sum\limits_{i=1}^{d_F} f_i^F$,
where $d_F \le n_{2j_F+1},f_1^F<\ldots<f_{d_F}^F$ in $G_0$ and
$E_{d_1}^F<E_{d_2}^F<E_{d_3}^F$ subintervals of
$\{1,\ldots,d_F\}$ such that
\begin{equation}\label{L.59}
\delta<\sum\limits_{q \in E_{d_i}^{F}}\ f_q^F(x_{n_{d_i}})
\mbox{~for all~} i=1,2,3.
\end{equation}
For every $F \in [L_1]^k$ we fix the functional $\phi_F$.
Throughout this proof, the functional $\phi_F$ will be called
the corresponding functional to $F$.\\
Let $L_1=\{l_1<\ldots<l_{d_1}<\ldots\}$
and $L_{1,d_1}=L_1 \backslash \{l_1<\ldots<l_{d_1}\}$.\\
Moreover if $F \in [L_{1,d_1}]^{k-d_1}$, we denote by
$\overline{F}$ the set $\overline{F}=\{l_1<\ldots<l_{d_1}\}\cup
F$.
\vspace{1.5mm}\\
\textbf{Claim:} There exists $C>0$ such that for every
$F \in [L_{1,d_1}]^{k-d_1}$ we have that
$w(\phi_{\overline{F}})\le C$.
\vspace{1.5mm}\\
Assume not. Then for $C=\frac{|\supp(x_{l_{d_1}})|}{\delta}$,
there exists $F \in [L_{1,d_1}]^{k-d_1}$ such that if
$\phi_{\overline{F}}$ is the corresponding functional to
$\overline{F}$, then from $(\ref{L.59})$ we get that
\[\delta<\sum\limits_{q \in E_{d_1}^{\overline{F}}}
f_q^{\overline{F}}(x_{l_{d_1}}).\]
On the other hand we have that
\[\sum\limits_{q \in E_{d_1}^{\overline{F}}}
f_q^{\overline{F}}(x_{l_{d_1}})\le |supp(x_{l_{d_1}})|\cdot
\sum\limits_{q \in E_{d_1}^{\overline{F}}} \frac{1}
{w(f_q^{\overline{F}})} \le \frac{|supp(x_{l_{d_1}})|}
{w(\phi_{\overline{F}})} <\frac{|supp(x_{l_{d_1}})|}
{C}=\delta \]
a contradiction.\\
Now applying twice Ramsey's theorem for the set
$[L_{1,d_1}]^{k-d_1}$ we obtain $L_2 \in [L_{1,d_1}],
j_0 \in \N$ and $r_1<r_2<r_3\in \{1,\ldots,n_{2j_0+1}\}$
such that
\begin{enumerate}
\item for every $F \in [L_2]^{k-d_1}$, it follows that
$w(\phi_{\overline{F}})=m_{2j_0+1}$ \item for every $F \in
[L_2]^{k-d_1}$ and every $i \in \{1,2,3\}$, we have that
$E_{d_i}^F=E_{d_i}$ \item $r_1 \in E_{d_1},r_2=\min E_{d_2},r_3
\in E_{d_3}$ and for every $F=\{n_{d_1+1}<\ldots<n_k\} \in
[L_2]^{k-d_1}$, then
\begin{enumerate}
\item [i.] $f_{r_1}^{\overline{F}}(x_{l_{d_1}})>
\frac{\delta}{n_{2j_0+1}},$
$f_{r_3}^{\overline{F}}(x_{n_{d_3}})>
\frac{\delta}{n_{2j_0+1}}$ and
\item [ii.] $f_{r_1}^{\overline{F}}<f_{r_2}^{\overline{F}}<
f_{r_3}^{\overline{F}}$.
\end{enumerate}
\end{enumerate}
We assume without loss of generality that $L_2=L_{1,d_1}$. Hence
$L_2=\{l_{d_1+1}<\ldots<l_{d_2-1}<\ldots\}$. Let $L_{2,d_2-1}=L_2
\backslash \{l_{d_1+1}<\ldots< l_{d_2-1}\}$. Moreover if $F \in
[L_{2,d_2-1}]^{k-(d_2-1)}$, we denote by $\overline{F}$ the set
$\overline{F}=\{l_1<\ldots<l_{d_1}<\ldots<l_{d_2-1}\}
\cup F$.\\
We consider the following sets
\[ A=\{F=\{n_{d_2}<\ldots<n_k\} \in [L_{2,d_2-1}]^{k-(d_2-1)}:
\#E_{n_{d_2}}^{\phi_{\overline{F}}} \ge 2\} \mbox{~and~}\]
\[B=\{F=\{n_{d_2}<\ldots<n_k\} \in [L_{2,d_2-1}]^{k-(d_2-1)}:
\#E_{n_{d_2}}^{\phi_{\overline{F}}}=1\}.\] It is obvious that
$[L_{2,d_2-1}]^{k-(d_2-1)}=A\cup B$. Hence from Ramsey's theorem
we may assume that $[L_{2,d_2-1}]^{k-(d_2-1)} \subset A$ or
$[L_{2,d_2-1}]^{k-(d_2-1)} \subset B$. We distinguish the
following cases.
\vspace{1.5mm}\\
\textbf{Case 1.} Let $[L_{2,d_2-1}]^{k-(d_2-1)} \subset A$.
\vspace{1.5mm}\\
Let $p \in \N$ such that $p>\frac{M \cdot n_{2j_0+1}}{\delta}$.\\
We consider vectors $y_1<\ldots<y_{p^2}$ such that
$y_i \in \{x_n:n \in L_{2,d_2-1}\}$ for all
$i=1,\ldots,p^2$ with $x_{l_{d_2-1}}<y_1$.
We set
\[y_i=x_{n_{{d_2}^{(i)}}},i=1,\ldots,p^2,\mbox{~where~}
n_{{d_2}^{(1)}}<\ldots<n_{{d_2}^{(p^2)}}
\mbox{~in~} L_{2,d_2-1} \mbox{~and~}
l_{d_2-1}<n_{{d_2}^{(1)}}.\]
We also fix vectors $x_{n_{d_2+1}}<\ldots<
x_{n_{d_3}}<\ldots<x_{n_k}$ such that
\[ n_{d_2+1}<\ldots<\ldots<n_{d_3}<\ldots<n_k \mbox{~in~}
L_{2,d_2-1} \mbox{~and~} y_{p^2}<x_{n_{d_2+1}}.\]
For every $i=1,\ldots,p^2$ we consider the following
subsets of $[L_{2,d_2-1}]^{k-(d_2-1)}$
\[F_i=\{n_{{d_2}^{(i)}}<n_{d_2+1}<\ldots<n_k\}.\]
For every $i=1,\ldots,p^2$ there exists
$f_{r_3}^{\overline{F_i}}\in G_0$ with even weight
such that $f_{r_3}^{\overline{F_i}}(x_{n_{d_3}})>
\frac{\delta}{n_{2j_0+1}}$.\\
Since $\#E_{n_{{d_2}^{(i)}}}^{\phi_{\overline{F_i}}} \ge 2$ for
all $i \in \{1,\ldots,p^2\}$, it follows that the functionals
$\{f_{r_2}^{\overline{F_i}}:i=1,\ldots,p^2\}$, are successive.
Hence the functionals
$\{f_{r_3}^{\overline{F_i}}:i=1,\ldots,p^2\}$,
have pairwise different weights.\\
We consider the functional
$f=\sum\limits_{i=1}^{p^2} \frac{1}{p}
f_{r_3}^{\overline{F_{l_i}}}$ which obviously belongs
to $G_0$. Hence \\
\[M<p\cdot \frac{\delta}{n_{2j_0+1}}\le f(x_{n_{d_3}})\le
\|x_{n_{d_3}}\|_{G_0} \le M\], a contradiction.
\vspace{1.5mm}\\
\textbf{Case 2.} Let $[L_{2,d_2-1}]^{k-(d_2-1)} \subset B$.
\vspace{1.5mm}\\
We consider the following sets:
\vspace{1mm}\\
$C=\{\{q_1<q_2\} \in [L_{2,d_2-1}]^2:\mbox{~if~}
F_1=\{q_1<n_{d_2+1}<\ldots<n_k\},
F_2=\{q_2<n_{d_2+1}<\ldots<n_k\}
\mbox{~subsets of~} [L_{2,d_2-1}]^{k-(d_2-1)},
\mbox{~then~} f_{r_2}^{\overline{F_1}}=
f_{r_2}^{\overline{F_2}}\}$
\vspace{1mm}
and $D=[L_{d_2-1}]^2 \backslash C$.\\
From Ramsey's Theorem we assume without loss of
generality that $[L_{d_2-1}]^2 \subset C$ or
$[L_{d_2-1}]^2 \subset D$.
We distinguish the following cases.
\vspace{1.5mm}\\
\textbf{Case 2.1.} Let $[L_{2,d_2-1}]^2 \subset D$.
\vspace{1.5mm}
In this case we derive to contradiction following the same
steps as in case 1.
\vspace{1mm}\\
\textbf{Case 2.2.} Let $[L_{2,d_2-1}]^2 \subset C$.
\vspace{1mm}\\
In this case we will prove that $\ell^1(\N)$ embeds isomorphically
in the space $\mathfrak{X}_{G_0}$, which contradicts to the fact
that $\mathfrak{X}_{G_0}$ is
reflexive.\\
We will prove that the sequence $(x_q)_{q \in L_{2,d_2-1}}$
is equivalent to the usual basis of $\ell^1(\N)$.\\
Let $n \in \N$ and $a_0,\ldots,a_n$ real numbers.
Since $L_{2,d_2-1}=\{l_{d_2+i},i=0,1,\ldots \}$
we get that $\|\sum\limits_{i=0}^{n} a_i x_{l_{d_2+i}}\|
_{G_0} \le M\sum\limits_{i=0}^{n} |a_i|.$ \\
We consider the following elements of
$[L_{2,d_2-1}]^{k-(d_2-1)}$
\[ F_q=\{q<n_{d_2+1}<\ldots<n_k\},q \in L_{2,d_2-1}.\]
From the fact that $[L_{2,d_2-1}]^2 \subset C$, we have that
$f_{r_2}^{\overline{F_{q_1}}}= f_{r_2}^{\overline{F_{q_2}}}=f$ for
all $q_1,q_2\in L_{2,d_2-1}$. Moreover, since
$\#E_{q}^{\phi_{\overline{F_q}}}=1$ for all $q \in L_{2,d_2-1}$,
we obtain that $E_{d_2}=E_{d_2}^{\overline{F_q}}=\{r_2\}$. Hence
from $(\ref{L.59})$ we have that
$f_{r_2}^{\overline{F_q}}(x_q)>\delta$ for all $q \in L_{2,d_2-1}$
and thus
\[ \|\sum\limits_{i=0}^{n} a_i x_{l_{d_2+i}}\|_{G_0}\ge
f(\sum\limits_{i=0}^{n} a_i x_{l_{d_2+i}})\ge
\sum\limits_{i=0}^{n} a_i
f(x_{l_{d_2+i}})\ge \delta
\sum\limits_{i=0}^{n} a_i \] and from the
1-unconditionality of the basis of the space
$\mathfrak{X}_{G_0}$ we obtain that
\[\|\sum\limits_{i=0}^{n} a_i x_{l_{d_2+i}}\|_{G_0}\ge
\delta \sum\limits_{i=0}^{n} |a_i|.\]

\vspace{3mm}
 The following is an immediate consequence of the
previous lemma.

\begin{Corollary}\label{L.10}
Let $k \in \N,\ \delta>0$ and $(x_n)_{n \in \N}$ be a block
sequence which is $\|\cdot\|_{W_0}$ bounded. Then for every $L \in
[\N]$, there exists $M \in [L]$ such that for every
$n_1<\ldots<n_k$ in $M$ the block finite sequence
$(x_{n_i})_{i=1}^{k}$ is not $\delta-$ separated.
\end{Corollary}

\begin{Lemma}\label{L.63}
Let $k,j_0 \in \N,\delta>0,(x_n)_{n=1}^k$ be a finite block
sequence and $\phi \in (W_0 \backslash G_0)$ of the form
$\phi=\frac{1}{m_{2j_0+1}}\sum\limits_{i=1}^{q} f_i$, where $q \le
n_{2j_0+1}$ and $(f_i)_{i=1}^{q}$ successive elements in $G_0$
such that
\begin{enumerate}
\item $(x_n)_{n=1}^k$ is not $(\phi,\delta)$ separated
and
\item $\delta<(\sum\limits_{i=1}^{q} f_i)(x_n)$ for all
$n=1,\ldots,k$.
\end{enumerate}
Then there exist at most four $d_1<d_2<d_3<d_4$ in
$\{1,\ldots,k\}$ such that $d_1=1,d_4=k$ and at most three
$i_1<i_2<i_3$ in $\{1,\ldots,q\}$ such that setting
\[I_s=\{n \in \{1,\ldots,k\}:d_s<n<d_{s+1}\} \mbox{~for~}
s=1,2,3\] it follows that
\[\phi(x_n)=\frac{1}{m_{2j_0+1}}f_{i_s}(x_n) \mbox{~for every~}
s=1,2,3 \mbox{~and every~} n \in I_s.\]
Hence \[\phi(x_1+\ldots+x_k)=\frac{1}{m_{2j_0+1}}\sum\limits_
{s=1}^{3} f_{i_s}(x_1+\ldots+x_k)+\sum\limits_{i=1}^{4}
\phi(x_{d_i}).\]
\end{Lemma}
\noindent{\bf{Proof:}}
Since $(x_n)_{n=1}^{k}$ is not $(\phi,\delta)$ separated,
we may assume without loss of generality that there exist
exactly two $r_1<r_2$ in $\{1,\ldots,k\}$ and
$E_{r_1}<E_{r_2}$ subintervals of $\{1,\ldots,q\}$ such that
$\delta<\sum\limits_{p \in E_{r_i}}f_p(x_{r_i})$ for all
$i=1,2$.\\
If there exists only one the proof is similar.\\
From the fact that $\delta<(\sum\limits_{i=1}^{q} f_i)(x_n)$ for
all $n=1,\ldots,k$, we obtain that
\[E_{r_1}\cap E_1^{\phi} \neq \emptyset \mbox{~and~}
E_{r_2}\cap E_k^{\phi} \neq \emptyset.\]
Setting
\[ A=\{d \in \{1,\ldots,k\}:\# E_d^{\phi} \ge 2\} \]
it is not hard to see that $\#A\le 4.$ We assume without loss of
generality that $\#A=4$. Then $A=\{d_1=1<r_1<r_2<d_4=k\}$. We set
$d_1=1,\ d_2=r_1,\ d_3=r_2$ and $d_4=k$. \\ Moreover setting
\[I_s=\{n \in \{1,\ldots,k\}:d_s<n<d_{s+1}\} \mbox{~for~}
s=1,2,3 \] we observe that the set
\[B=\{d \in \{1,\ldots,k\}: \# E_d^{\phi}=1\}=
\bigcup_{s=1}^{3} I_s.\]
Therefore for every $s=1,2,3$ there exists $i_s \in
\{1,\ldots,q\}$ with $i_1<i_2<i_3$ such that
the conclusion of lemma is satisfied.

\begin{Lemma}\label{L.26}
Let $k \in \N$, $\delta>0$ with $\sqrt{k}>\frac{4}{\delta}$ and
$j_0 \in \N$. Let also $(x_n)_{n=1}^{k}$ be a finite block
sequence with $\|x_n\|_{W_0} \le 1$ for all $n=1,\ldots,k$ and
$\phi \in(W_0 \backslash G_0)$ with $w(\phi)=m_{2j_0+1}$ such that
\begin{enumerate}
\item $2\delta<\|\frac{x_1+\ldots+x_k}{k}\|_{G_0}$ and
\item $(x_n)_{n=1}^{k}$ is not $(\phi,\delta)$ separated.
\end{enumerate}
Then
\[\phi(\frac{x_1+\ldots+x_k}{k})\le
(\frac{4}{m_{2j_0+1}}+\frac{1}{2\sqrt{k}})
\|\frac{x_1+\ldots+x_k}{k}\|_{G_0}.\]
\end{Lemma}
\noindent{\bf{Proof:}} Let $\phi$ be of the form
$\phi=\frac{1}{m_{2j_0+1}}\sum\limits_{i=1}^{q} f_i$, where $q \le
n_{2j_0+1}$ and $(f_i)_{i=1}^{q}$ successive elements in $G_0$.
Let
\[D=\{n \in \{1,\ldots,k\}:\delta<
(\sum\limits_{i=1}^{q} f_i)(x_n) \}\]
and
\[D_0=\{1,\ldots,k\} \backslash D.\]
Since $(x_n)_{n=1}^{k}$ is not $(\phi,\delta)$ separated, we get
that $(x_n)_{n \in D}$ is also not $(\phi,\delta)$ separated.
Next, using Lemma $\ref{L.63}$ we will estimate the real number
$\phi(\sum\limits_{i=1}^{k} x_i)$. We have that
\begin{eqnarray*}
\phi(\sum\limits_{i=1}^{k} x_i) & & =
\phi(\sum\limits_{i \in D_0} x_i)+
\phi(\sum\limits_{i\in D} x_i)\le
\frac{\delta \cdot k}{m_{2j_0+1}}+
\frac{1}{m_{2j_0+1}}\sum\limits_{s=1}^{3}
f_{i_s}(x_1+\ldots+x_k)+
\sum\limits_{i=1}^{4}\phi(x_{d_i})
\\ & & \le \frac{1}{m_{2j_0+1}}
\| \sum\limits_{i=1}^{k} x_i\|_{G_0}+
\frac{3}{m_{2j_0+1}}\|\sum\limits_{i=1}^{k} x_i\|_{G_0}
+4 \le
\frac{4}{m_{2j_0+1}}\|\sum\limits_{i=1}^{k} x_i\|_{G_0}+4
\end{eqnarray*}
Hence
\begin{eqnarray*}
\phi(\frac{x_1+\ldots+x_k}{k}) & &=
\frac{4}{m_{2j_0+1}}\|\frac{x_1+\ldots+x_k}{k}\|_{G_0}
+\frac{4}{\sqrt{k}} \cdot \frac{1}{\sqrt{k}} \le
\frac{4}{m_{2j_0+1}}\|\frac{x_1+\ldots+x_k}{k}\|_{G_0}
+\frac{\delta}{\sqrt{k}} \\ & &
\le \frac{4}{m_{2j_0+1}}\|\frac{x_1+\ldots+x_k}{k}\|_{G_0}
+\frac{1}{2\sqrt{k}}\|\frac{x_1+\ldots+x_k}{k}\|_{G_0}
\\ & &
=(\frac{4}{m_{2j_0+1}}+\frac{1}{2\sqrt{k}})
\|\frac{x_1+\ldots+x_k}{k}\|_{G_0}.
\end{eqnarray*}
\vspace{2mm}\\

The previous result yields immediately the following.

\begin{Corollary}\label{L.60}
Let $k \in \N$ and $\delta>0$ with
$\sqrt{k}>\frac{4}{\delta}$.
Let also $(x_n)_{n=1}^{k}$ be a finite block
sequence with $\|x_n\|_{W_0} \le 1$ for all
$n=1,\ldots,k$ and
$\phi \in(W_0 \backslash G_0)$ such that
\[2\delta<\|\frac{x_1+\ldots+x_k}{k}\|_{G_0}<
\phi(\frac{x_1+\ldots+x_k}{k}).\]
Then $(x_n)_{n=1}^{k}$ is $(\phi,\delta)$ separated.
\end{Corollary}

\begin{Proposition}\label{L.24}
For every $\epsilon>0$ there exists $n \in \N$
such that for every $k \in \N$ with $k>n$
and every block sequence $(x_n)_{n \in \N}$ with
\[0<\epsilon<\|x_n\|_{G_0}\le \|x_n\|_{W_0}\le
1 \mbox{~for all~} n \in \N\]
there exists an $L\in [\N]$ such that
\\ for every $n_1<\ldots<n_k$ in $L$ and every
$\phi\in (W_0\backslash G_0)$ it follows that
\[\phi(\frac{x_{n_1}+\ldots+x_{n_k}}{k})\le
\|\frac{x_{n_1}+\ldots+x_{n_k}}{k}\|_{G_0}.\]
Hence
\[\|\frac{x_{n_1}+\ldots+x_{n_k}}{k}\|_{W_0}=
\|\frac{x_{n_1}+\ldots+x_{n_k}}{k}\|_{G_0}.\]
\end{Proposition}
\noindent{\bf{Proof:}}
Assume that the conclusion of the proposition fails.
Then there exists $\epsilon>0$ such that for
every $n \in \N$ there exist $k>n$ and a block sequence
$(x_n)_{n \in \N}$ with
\[\epsilon<\|x_n\|_{G_0}\le \|x_n\|_{W_0}\le 1
\mbox{~for all~} n \in \N \] such that for every $L\in [\N]$,
there exist $n_1<\ldots<n_k$ in $L$ and $\phi\in (W_0 \backslash
G_0)$ with
\[\|\frac{x_{n_1}+\ldots+x_{n_k}}{k}\|_{G_0}<
\phi(\frac{x_{n_1}+\ldots+x_{n_k}}{k}).\]
There exist $k,j_k$ in $\N$ such that
\begin{enumerate}
\item  $m_{2j_k}<\sqrt{k}<k<n_{2j_k}$
\item  $\frac{8}{\epsilon}<\frac{\sqrt{k}}{m_{2j_k}}$.
\end{enumerate}
and a block sequence $(x_n)_{n \in \N}$ satisfying the above
properties. Hence from Ramsey's theorem there exists $L \in [\N]$
such that for every $n_1<\ldots<n_k$ in $L$ there exists $\phi\in
(W_0 \backslash G_0)$ with
\begin{equation}\label{L.65}
\|\frac{x_{n_1}+\ldots+x_{n_k}}{k}\|_{G_0}<
\phi(\frac{x_{n_1}+\ldots+x_{n_k}}{k}).
\end{equation}
We observe from the choice of $k \in \N$, that for every
$n_1<\ldots<n_k$ in $L$ we obtain
\begin{equation}\label{L.66}
\frac{\epsilon}{m_{2j_k}}<\|\frac{x_{n_1}+
\ldots+x_{n_k}}{k}\|_{G_0}.
\end{equation}
From $(\ref{L.65}), (\ref{L.66})$ and Corollary \ref{L.60} we get
that for every $F \in [L]^k$, the sequence $(x_n)_{n \in F}$ is
$\frac{\epsilon}{m_{2j_k}}-$ separated. Thus from Corollary
\ref{L.10}, we derive to contradiction.

\begin{Proposition}\label{L.25}
Let $\epsilon>0$ and
$(x_n)_{n \in \N}$ be a block sequence such that
\[\epsilon<\|x_n\|_{G_0}\le \|x_n\|_{W_0}\le 1 \mbox{~for all~}
n \in \N.\]
Then for every $n \in \N$ there exists $k \in \N$ with $k>n$
and there exist $y_1<\ldots<y_k$ in $<x_n:n \in \N>$ such that
\begin{enumerate}
\item $\|y_i\|_{W_0} \le 1$ for all $i=1,\ldots,k$ and
\item if $y=\frac{1}{k}(y_1+\ldots+y_k)$ then $\frac{1}{2}
<\|y\|_{G_0}=\|y\|_{W_0}$.
\end{enumerate}
\end{Proposition}
\noindent{\bf{Proof:}} Assume that the conclusion fails. Then
there exists $n \in \N$ such that for every $k \in \N$ with $k>n$
and every $y_1<\ldots<y_k$ in $<x_n:n \in \N>$ with
$\|y_i\|_{W_0}\le 1$ for all $i=1,\ldots,k$ it follows that
$\|\frac{1}{k}(y_1+\ldots+y_k)\|_{G_0}\le\frac{1}{2}$ or
$\|\frac{1}{k}(y_1+\ldots+y_k)\|_{G_0}<
\|\frac{1}{k}(y_1+\ldots+y_k)\|_{W_0}.$ \\
It is obvious that there exists $k_0 \in \N$ with $k_0>n$
such that for every $k \in \N$ with $k>k_0$
the conclusion of Proposition $\ref{L.24}$ is valid.\\
We choose $j,s \in \N$ such that
$m_{2j}\le \epsilon \cdot 2^{s-1}<2^s\le k_0^s \le n_{2j}$.\\
We set $z_n^{(1)}=x_n,n \in \N$. Then from Proposition
$\ref{L.24}$ there exists an $L_1 \in [\N]$ such that for every
$n_1<\ldots<n_{k_0}$ in $L_1$ it follows that
\begin{equation}\label{L.17}
\|\frac{1}{k_0}(z_{n_1}^{(1)}+\ldots+z_{n_{k_0}}^{(1)})\|_{G_0}=
\|\frac{1}{k_0}(z_{n_1}^{(1)}+\ldots+z_{n_{k_0}}^{(1)})\|_{W_0}.
\end{equation}
Let $L_1=\{l_1^{(1)}<\ldots<l_i^{(1)}<l_{i+1}^{(1)}<\ldots\}$.
We set
\[ w_{1,n}=\sum\limits_{i=(n-1)k_0+1}^{nk_0} x_{{l_i}^{(1)}},
n \in \N. \] Let $n \in \N$. Since $\|z_n^{(1)}\|_{W_0}\le 1$ and
$(\ref{L.17})$ holds, we obtain that
$\|\frac{1}{k_0}w_{1,n}\|_{G_0}=
\|\frac{1}{k_0}w_{1,n}\|_{W_0}\le\frac{1}{2}$. Hence
$\|\frac{2}{k_0}w_{1,n}\|_{W_0}\le 1$.
\\ We set
$z_n^{(2)}=\frac{2}{k_0}w_{1,n},n \in \N.$
Then from Proposition $\ref{L.24}$ there exists an
$L_2 \in [\N]$ such that for every $n_1<\ldots<n_{k_0}$ in
$L_2$ it follows that
\begin{equation}\label{L.18}
\|\frac{1}{k_0}(z_{n_1}^{(2)}+\ldots+z_{n_{k_0}}^{(2)})\|_{G_0}=
\|\frac{1}{k_0}(z_{n_1}^{(2)}+\ldots+z_{n_{k_0}}^{(2)}\|_{W_0}.
\end{equation}
Let $L_2=\{l_1^{(2)}<\ldots<l_i^{(2)}<l_{i+1}^{(2)}<\ldots\}$.
We set
\[w_{2,n}=\sum\limits_{i=(n-1)k_0+1}^{nk_0} w_{1,{l_i}^{(2)}},
n \in \N.\]
It is obvious that for every $n \in \N$ the vector $w_{2,n}$
consists of $k_0^2$ blocks of the sequence
$(x_n)_{n \in L_1}$.\\
Moreover, since for a $n \in \N$ we have that
$\|\frac{2}{k_0}w_{1,n}\|_{W_0}\le 1$, we get from $(\ref{L.18})$
that $\|\frac{2}{k_0^2}w_{2,n}\|_{G_0}=
\|\frac{2}{k_0^2}w_{2,n}\|_{W_0}\le \frac{1}{2}$. Hence
$\|\frac{2^2}{k_0^2}w_{2,n}\|_{W_0}\le 1$ for all $n \in \N$.\\In
this way we inductively construct
\begin{enumerate}
\item a finite sequence $(L_r)_{r=1}^{s}$ of infinite subsets of
$\N$, where
\\ $L_r=\{l_1^{(r)}<\ldots<l_i^{(r)}<l_{i+1}^{(r)}<\ldots\},
r=1,\ldots,k_0^s$ and
\item for every $r=1,\ldots,s$ a sequence
$(w_{r,n})_{n \in \N}$ of the form
$w_{r,n}=\sum\limits_{i=(n-1)k_0+1}^{nk_0} w_{r-1,{l_i}^{(r)}},
n \in \N,r=2,\ldots,s$ where
$w_{1,n}=\sum\limits_{i=(n-1)k+1}^{nk_0} x_{{l_i}^{(1)}},
n \in \N$ such that \\
$\|w_{r,n}\|_{G_0}=\|w_{r,n}\|_{W_0} \le\frac{k_0^r}{2^r}$
for $r=1,\ldots,s.$
\end{enumerate}
We can see that if $r \in \{1,\ldots,s\}$ then for every $n \in
\N$ the vector $w_{r,n}$ consists of $k_0^r$ block vectors of the
sequence $(x_n)_{n \in \L_1}$.\\ Let $n_0 \in \N$ and
$w_{s,n_0}=\sum\limits_{i=1}^{k_0^s} x_{l_{q_i}}^{(1)}$. For every
$i=1,\ldots,k_0^s$ there exists $g_i \in G_0$ with
$ran(g_i)\subset ran(x_{l_{q_i}}^{(1)})$ such that
$g_i(x_{l_{q_i}}^{(1)})>\frac{\epsilon}{2}$. We consider the
functional $g=\frac{1}{m_{2j}}\sum\limits_{i=1}^{k_0^s} g_i$ which
belongs to $G_0$. Then
\[\frac{\epsilon\cdot k_0^s}{2m_{2j}}<g(w_{s,n_0})\le
\frac{k_0^s}{2^s}\] which yields that $2^{s-1}<m_{2j}$,
a contradiction.

\section{\textbf{Rapidly increasing sequences in
$\mathfrak{X}_{G_0}$ and in $\mathfrak{X}_0$}}

We begin with the definition of the Rapidly increasing
Sequences(RIS) in $\mathfrak{X}_{G_0}$ and in $\mathfrak{X}_0$ and
the definition of $M-\ell_k^1$ averages in $\mathfrak{X}_0$.
\vspace{3mm}

\begin{Definition}(\textbf{Rapidly increasing sequences in
$\mathfrak{X}_{G_0},\mathfrak{X}_0$}) Let $(x_n)_{n\in \N}$ be a
block sequence and $C,\epsilon$
positive numbers. This sequence will be called \\
\textbf{(i.)} $(C,\epsilon)$ $RIS$ in $\mathfrak{X}_{G_0}$
if the following hold:
\begin{enumerate}
\item[1.] $\|x_n\|_{G_0} \leq C$ for all $n\in\N$.
\item[2.] There exists a strictly increasing sequence of
natural numbers $(j_n)_{n\in \N}$ such that \\
$\frac{|\supp x_n|}{m_{j_{n+1}}}<\epsilon$, for all $n\in\N$.
\item[3.] For every $n\in\N$ and every $f\in G_0$ with
$w(f)=m_{2i}<m_{j_n}$ we have that $|f(x_n)|\leq\frac{C}{m_{2i}}$.
\end{enumerate}
\textbf{(ii.)} $(C,\epsilon)$ $RIS$ in $\mathfrak{X}_0$ if the
following hold:
\begin{enumerate}
\item[1.] $\|x_n\|_{W_0} \leq C$ for all $n\in\N$.
\item[2.] There
exists a strictly increasing sequence of natural numbers
$(j_n)_{n\in \N}$ such that \\
$\frac{|\supp x_n|}{m_{j_{n+1}}}<\epsilon$, for all $n\in\N$.
\item[3.] For every $n\in\N$ and  $f\in W_0$ with
$w(f)=m_i<m_{j_n}$ we have that $|f(x_n)|\leq\frac{C}{m_i}$.
\end{enumerate}
Every sequence $(j_n)_{n \in \N}$ of natural numbers like in this
definition is said to be the associated sequence of the $RIS$
$(x_n)_{n \in \N}$.
\end{Definition}
\vspace{2mm}

Next, we define the $M-\ell_k^1$ averages in $\mathfrak{X}_0$.

\begin{Definition}(\textbf{$\ell_k^1-$ averages in
$\mathfrak{X}_0$}) Let $k\in\N,\ M>0$ and $(e_n)_{n \in \N}$ the
Schauder basis of $\mathfrak{X}_0$.\\
A vector $x\in <e_n:\;n \in \N>$ is said to be a
$M-\ell_k^1$ average in $\mathfrak{X}_0$ if:
\begin{enumerate}
\item[1.] $\|x\|_{G_0}=\|x\|_{W_0}>\frac{1}{2}.$
\item[2.] There exists
$x_1<\ldots<x_k$ in $<e_n:\;n \in \N>$ with
$\|x_i\|_{W_0}\leq M$ for all $i=1,...,k$ such that
$x=\frac{1}{k}\sum_{i=1}^k x_i.$
\end{enumerate}
\end{Definition}
\vspace{4mm}

In the following Lemma, according to Proposition $\ref{L.25}$, it
is easily checked that in every block subspace $<x_n,n \in \N>$,
where
\[0<\epsilon<\|x_n\|_{G_0}\le \|x_n\|_{W_0}\le 1
\mbox{~for all~} n \in \N\]
there exists a block sequence of
$(x_n)_{n \in \N}$ of $l^1-$ averages with increasing lengths.
A consequence of this result is the existence of a
$(3,\delta)$ RIS in $\mathfrak{X}_0$ for a fixed $\delta>0$.
The proof of the existence of the RIS follows the lines of
the proof of Proposition II.25 in \cite{ATO}.

\begin{Lemma}\label{L.1}
Let $\epsilon,\ \delta>0$ and $(x_n)_{n\in \N}$ be a block
sequence such that
\[\epsilon<\|x_n\|_{G_0}\le \|x_n\|_{W_0}\le 1
\mbox{~for all~} n \in \N.\]
Then there exist a block sequence
$(y_n)_{n\in \N}$ of $(x_n)_{n\in \N}$ and a subsequence
$(k_n)_{n \in \N}$ of natural numbers such that
\begin{enumerate}
\item[i.] For every $n\in \N$ the vector $y_n$ is a
$1-\ell_{k_n}^1$ average in  $\mathfrak{X}_0$.
\item[ii.] The sequence $(y_n)_{n\in \N}$ is $(3,\delta)$
RIS in $\mathfrak{X}_0$.
\end{enumerate}
\end{Lemma}
\noindent{\bf{Proof:}} The first assertion is an immediate
consequence of Proposition $\ref{L.25}$. The proof of the second
follows the lines of the proof of Proposition II.25 in \cite{ATO}.

\section{\textbf{The existence of exact pairs in
$\mathfrak{X}_0$}}

The aim of this section is the existence of exact pairs and
dependent sequences in $\mathfrak{X}_0$. In order to achieve this
we use the Basic Inequality for the space $\mathfrak{X}_0$. As
usually we will use an auxiliary space $F_{j_0}$, $j_0 \in \N$ and
our approach follows similar steps as in \cite{AMP}. \vspace{2mm}

\begin{Definition}\label{L.67}
Let $j_0\in \N$ with $j_0>1$. We denote by $F_{j_0}$ the minimal
subset of $c_{00}(\N)$ such that:
\begin{enumerate}
\item $C_{j_0}=\{\sum\limits_{i \in F}
\epsilon_ie_i^*:|\epsilon_i|=1, \#F \le n_{j_0-1}\} \subset
F_{j_0}$. \item $F_{j_0}$ is closed under the operation
$(\mathcal{A}_{2n_{2j}},\frac{1}{m_{2j}})$ for every $j\in \N$.
\item For every $A,B$ nonempty finite subsets of $\N$ with $A\cap
B=\emptyset$, for every $(a_i)_{i \in A \cup B}$ finite sequence
of real numbers such that $\sum\limits_{i \in A \cup B} a_i^2 \le
1$, for every $(f_i)_{i \in A}$, where $f_i \in F_{j_0}$ is a
result of a $(\mathcal{A}_{2n_{j_i}}, \frac{1}{m_{j_i}})$
operation with $w(f_{i_1})\neq w(f_{i_2})$ for all $i_1\neq i_2
\in A$, and for every finite sequence of natural numbers $(t_i)_{i
\in B}$, where $t_i\neq t_j$ for all $i\neq j \in B$, it follows
that $(\sum\limits_{i \in A} a_if_i+ \sum\limits_{i \in B}
a_ie_{t_i}) \in F_{j_0}$.
\end{enumerate}
\vspace{3mm} We notice that for an $f \in F_{j_0}$ we say that $f$
has weight $m_{2j}$ or that is a result of a
$(\mathcal{A}_{2n_{2j}}, \frac{1}{m_{2j}})$ operation and we write
$w(f)=m_{2j}$ if and only if there exists $d \in \N$ with $d\le
2n_{2j}$  and $f_1<\ldots<f_d \in F{j_0}$ such that
$f=\frac{1}{m_{2j}}\sum_{i=1}^d f_i$.\\
We also define the auxiliary space $F_{j_0}^{\prime}$ to be
the following subset of $c_{00}(\N)$.
\[F_{j_0}^{\prime}=F_{j_0} \cup
\{\frac{1}{m_{2j+1}}\sum_{i=1}^d f_i:j \in \N,d \le 2n_{2j+1},
f_i \in F_{j_0}\}.\]
The weight of a functional $f$ belonging to $F_{j_0}^{\prime}
\backslash F_{j_0}$ is defined in a similar way as for the
functionals of $F_{j_0}$.
\end{Definition}

We begin with the estimations of a functional belonging to the
auxiliary space, acting in averages of the basis of length $n_j$.
For a proof of the following Lemma we refer to Lemma 8.10 in
\cite{AMP}.

\begin{Lemma}\label{L.7}
There exists $M>0$ such that if $j_0\in \N$ then
\begin{enumerate}
\item[1.] if $f \in F_{j_0}^{\prime}$ with
$w(f)=m_i,i\in \N$ and
$k_1<\ldots<k_{n_{j_0}}$ natural numbers then
\[|f(\frac{1}{n_{j_0}}\sum\limits_{r=1}^{n_{j_0}}
e_{k_r})| \le \begin{cases} \frac{4}{m_{i}m_{j_0}},~
\mbox{if~} i<j_0 \\
\frac{M}{m_{i}},~ \mbox{if~} i\ge j_0
\end{cases}\]
\item[2.] if $f \in F_{2j_0+1}^{\prime}$ with
$w(f)=m_i,i \neq 2j_0+1$ and
$k_1<\ldots<k_{n_{2j_0+1}}$ natural numbers then
\[|f(\frac{1}{n_{2j_0+1}}\sum\limits_{r=1}^{n_{2j_0+1}}
e_{k_r})| \le \begin{cases}
\frac{4}{m_im_{2j_0+1}^2},~ \mbox{if~} i<2j_0+1 \\
\frac{M}{m_i},~ \mbox{if~} i\ge 2j_0+1
\end{cases}\]
\end{enumerate}
\end{Lemma}

\vspace{2mm}
The proof of the following proposition (basic inequality
in $\mathfrak{X}_0$) follows the lines of the proof of
Proposition 9.3 in \cite{AMP}.

\begin{Proposition}  ({\bf basic inequality in $\mathfrak{X}_0$})
Let $(x_k)_{k\in \N}$ be a $(C,\epsilon)$ RIS in $\mathfrak{X}_0$
where $C>0,\ \epsilon>0$ and $j_0 \in \N$. Let also
$(\lambda_n)_{n\in \N}$ be a sequence of real numbers.\\
Then for
every $f\in W_0$ and every finite interval $I$ of $\N$ there
exist a functional $g\in F_{j_0}^{\prime}$ and a real number
$\epsilon_f\le \epsilon$
such that \[|f(\sum\limits_{k \in I} \lambda_kx_k)| \le
C(g(\sum\limits_{k \in I} |\lambda_k|e_k)+
\epsilon_f\cdot \sum\limits_{k \in I} |\lambda_k|). \]
If we assume that $w(f)=m_{2j}$ then $g=0$ or $g\in
\{\pm e_n^*:n \in \N\}$ or $w(g)=w(f)$ and $\epsilon_f
\le \frac{\epsilon}{w(f)}$.
\end{Proposition}

\begin{Remark}
The Basic Inequality in $\mathfrak{X}_{G_0}$ is analogous to that
in $\mathfrak{X}_0$. Specifically, for a block $(C,\epsilon)$ RIS
 in $\mathfrak{X}_{G_0}$ and a functional $f\in G_0$, there
exist a functional $g\in F_{j_0}$ and a real number $\epsilon_f\le
\epsilon$ satisfying the same properties as for the basic
inequality in $\mathfrak{X}_0$.
\end{Remark}
\vspace{2mm}

The proofs of the following two propositions (Proposition
$\ref{L.8}$ and $\ref{L.4}$) are based in Basic Inequality in
$\mathfrak{X}_0$ and Lemma \ref{L.7}. For a proof we refer to
Proposition 9.4 in \cite{AMP}.

\begin{Proposition}\label{L.8}
Let $j_0 \in \N$ and $(x_n)_{n\in \N}$ be a $(C,\epsilon)$ RIS in
$\mathfrak{X}_0$ with $(j_n)_{n \in \N}$ its associated sequence
such that $C>0$, $0<\epsilon<\frac{2}{m_{j_0}^2}$ and $j_1>j_0$.
Let also $k_1<\ldots<k_{n_{j_0}}$ natural numbers and $f \in W_0$
with $w(f)=m_i$. Then
\begin{enumerate}
\item[i.] $|f(\frac{x_{k_1}+\ldots+x_{k_{n_{j_0}}}}{n_{j_0}})|
\le
\begin{cases}
\frac{5C}{m_im_{j_0}}, \mbox{~if~} i<j_0 \\
\frac{MC}{m_{i}}+\frac{2C}{m_{j_0}^2}, \mbox{~if~} i\ge j_0 \\
\end{cases}$
\item[ii.] $\|\frac{x_{k_1}+\ldots+x_{k_{n_{j_0}}}}{n_{j_0}}\|_
{W_0} \le \frac{3C}{m_{j_0}}.$
\end{enumerate}
$M$ is the positive number appearing in Lemma \ref{L.7}.
\end{Proposition}

\begin{Remark}
For a $(C,\epsilon)$ RIS in $\mathfrak{X}_{G_0}$ we have the same
estimations for every functional $f$ belonging to $G_0$ and the
$\|\cdot\|_{G_0}$ norm as in Proposition $\ref{L.8}$.
\end{Remark}

\begin{Proposition}\label{L.4}
Let $j_0 \in \N$ and $(x_n)_{n\in \N}$ be a $(C,\epsilon)$ RIS in
$\mathfrak{X}_{G_0}$ with $(j_n)_{n \in \N}$ its associated
sequence such that $C>0$, $0<\epsilon<\frac{2}{m_{2j_0+1}^2}$ and
$j_1>2j_0+1$. Let also $k_1<\ldots<k_{n_{2j_0+1}}$ natural numbers
and $|a_i|\le 1$ for all $i=1,\ldots,n_{2j_0+1}$. Then
\[\|\frac{a_1x_{k_1}+\ldots+a_{n_{2j_0+1}}
x_{k_{n_{2j_0+1}}}}{n_{2j_0+1}}\|_{G_0} \le
\frac{3C}{m_{2j_0+1}^2}.\]
\end{Proposition}

\vspace{2mm} Next, we define the exact pairs in $\mathfrak{X}_0$.

\begin{Definition}(\textbf{exact pairs in $\mathfrak{X}_0$})
Let $x\in \mathfrak{X}_0$ with finite support and $f\in G_0$. The
pair $(x,f)$ is called a $(C,2j,\theta)$ exact pair in
$\mathfrak{X}_0$, where $C \geq 1,j\in\N$ and $\theta\ge 0$ if the
following hold:
\begin{enumerate}
\item[(i)] $\frac{1}{2}\le\|x\|_{G_0}=\|x\|_{W_0} \leq C$
\item[(ii)] $w(f)=m_{2j}$ \item[(iii)] $f(x)=\theta$ \item[(iv)]
If $g \in W_0$ with $w(g)=m_i$, then $|g(x)|\leq\frac{C}{m_i}$ if
$i<2j$ and $|g(x)|\le \frac{C}{m_{2j}}$, if $i>2j$.
\end{enumerate}
\end{Definition}

\vspace{4mm} In the following Proposition, using Proposition
$\ref{L.24}$, Lemma $\ref{L.1}$ and Basic Inequality in
$\mathfrak{X}_0$, we get that in every block subspace $<x_n,n \in
\N>$, where
\[0<\epsilon<\|x_n\|_{G_0}\le \|x_n\|_{W_0}\le 1
\mbox{~for all~} n \in \N\]
there exist an exact pair in $\mathfrak{X}_0$.
The proof follows the lines of the proof of
Proposition II.32 in \cite{ATO} or Proposition
10.2 in \cite{AMP}.

\begin{Proposition}\label{L.20}
Let $j \in \N,\epsilon>0$ and $(x_n)_{n\in \N}$ be
a block sequence such that
\[\epsilon<\|x_n\|_{G_0}\le \|x_n\|_{W_0}\le 1
\mbox{~for all~} n \in \N.\]
Then there exists a pair $(x,f)$ with
$x \in <x_n,n \in \N>$ which is $(15,2j,\frac{1}{2})$
exact pair in $\mathfrak{X}_0$.
\end{Proposition}
\noindent{\bf{Proof:}} From Lemma $\ref{L.1}$, there exist a block
sequence $(y_n)_{n\in \N}$ of $(x_n)_{n\in \N}$ and a subsequence
$(k_n)_{n \in \N}$ of natural numbers such that
\begin{enumerate}
\item for every $n\in \N$ the vector $y_n$ is a $1-\ell_{k_n}^1$
average in  $\mathfrak{X}_0$ and \item $(y_n)_{n \in \N}$ is
$(3,\delta)$ RIS in $\mathfrak{X}_0$, where
$0<\delta<\frac{2}{m_{2j}^2}$ and the first term of the associated
sequence is bigger than $2j$.
\end{enumerate}
Therefore if $k_1<\ldots<k_{n_{2j}}$ natural numbers,
then there exist
$f_1<\ldots<f_{n_{2j}}$ in $G_0$ such that
$f_i(y_{k_i})=\frac{1}{2}$ and $ran(f_i) \subset
ran(y_{k_i})$ for all $i=1,\ldots,n_{2j}$.
We set
\[x=\frac{m_{2j}}{n_{2j}}\sum\limits_{i=1}^{n_{2j}}y_{k_i}
\mbox{~and~} f=\frac{1}{m_{2j}} \sum\limits_{i=1}^{n_{2j}}f_i.\]
Then it easily checked, using Proposition $\ref{L.8}$, that the
pair $(x,f)$ is a $(15,2j,\frac{1}{2})$ exact pair in
$\mathfrak{X}_0$.

\vspace{2mm}
We define the dependent sequences in $\mathfrak{X}_0$.

\begin{Definition}(\textbf{dependent sequences in
$\mathfrak{X}_0$})
Let $j\in \N,\theta \ge 0$ and $C\ge 1$.
A finite sequence of pairs
$(x_i,f_i)_{i=1}^{n_{2j+1}}$ with $x_i \in \mathfrak{X}_0$
for all $i \in \{1,\ldots,n_{2j+1}\}$ is called a
$(C,2j+1,\theta)$-dependent sequence in $\mathfrak{X}_0$ if the
following hold:
\begin{enumerate}
\item[1.] $(f_i)_{i=1}^{n_{2j+1}}$ is a $\sigma-n_{2j+1}$ special
sequence with $w(f_i)=m_{2j_i}$ where
$2j_i=\sigma(f_1,\ldots,f_{i-1}),i\in \{2,\ldots,n_{2j+1}\}$ and
$2j_1 \in \{2j:\;j \in \Omega_1\}$ \item[2.] each pair $(x_i,f_i)$
is a $(C,2j_i,\theta)$ exact pair in $\mathfrak{X}_0$ \item[3.]
$\ran(f_i) \cup \ran(x_i)<\ran(f_{i+1}) \cup \ran(x_{i+1})$ for
all $i=1,\ldots,n_{2j+1}-1$.
\end{enumerate}
\end{Definition}

\vspace{4mm} In the following Proposition, using the existence of
exact pairs in every block subspace $<x_n,n \in \N>$ of
$\mathfrak{X}_0$, with
\[0<\epsilon<\|x_n\|_{G_0}\le \|x_n\|_{W_0}\le 1
\mbox{~for all~} n \in \N\]
we construct a dependent sequence in $\mathfrak{X}_0$.

\begin{Proposition}\label{L.9}
Let $j \in \N,\epsilon>0$ and $(x_n)_{n\in \N}$ be
a block sequence such that
\[\epsilon<\|x_n\|_{G_0}\le \|x_n\|_{W_0}\le 1
\mbox{~for all~} n \in \N.\] Then there exists a finite sequence
of pairs $(x_i,f_i)_{i=1}^{n_{2j+1}}$ with $x_i \in <x_n,n \in
\N>$ for all $i=1,\ldots,n_{2j+1}$, which is a
$(15,2j+1,\frac{1}{2})$-dependent sequence in $\mathfrak{X}_0$
with $\ran(f_i)=\ran(x_i)$ for all $i=1,\ldots,n_{2j+1}$.
\end{Proposition}
\noindent{\bf{Proof:}}
It follows easily from an inductive application of
Proposition $\ref{L.20}$.

\begin{Proposition}\label{L.2}
Let $j_0 \in \N,\ \theta\ge 0,\ C\ge 1$ and
$(x_i,f_i)_{i=1}^{n_{2j_0+1}}$ be a $(C,2j_0+1,\theta)$-dependent
sequence in $\mathfrak{X}_0$ with $maxsupp(f_i)\ge \#supp(x_i)$
for all $i=1,\ldots,n_{2j+1}$. Then
\begin{enumerate}
\item if $\theta=\frac{1}{2}$, it holds
\[\|\frac{1}{n_{2j_0+1}} \sum\limits_{i=1}^{n_{2j_0+1}}
(-1)^{i}x_i\|_{W_0} \le \frac{8C}{m_{2j_0+1}^2}.\] \item if
$\theta=0$, it holds
\[\|\frac{1}{n_{2j_0+1}} \sum\limits_{i=1}^{n_{2j_0+1}}
x_i\|_{W_0} \le \frac{8C}{m_{2j_0+1}^2}.\]
\end{enumerate}
\end{Proposition}
\noindent{\bf{Proof:}}
The proof follows the lines of the proof of Proposition
III.6 in \cite{ATO}.

\section{\textbf{Properties of $\mathfrak{X}_0$}}

In this section we prove that the Banach space $\mathfrak{X}_0$ is
$c_0-$ saturated, i.e. every closed, infinite dimensional subspace
contains an isomorphic copy of $c_0$. Also it is proved that if
$Y,Z$ are closed, infinite dimensional subspaces of
$\mathfrak{X}_0$ such that the direct sum $Y\bigoplus Z$ is a
closed subspace, then one at least of the subspaces is embedded
isomorphically into $c_0$ and the other contains an isomorph of
$\mathfrak{X}_0$. Our proof uses a deep result due to Kalton that
characterizes the subspaces of $c_0$. For the convenience, by a
subspace we always mean a closed and infinite dimensional one.

\vspace{2mm}

\begin{Proposition}\label{L.92}
The identity operator
$id:\mathfrak{X}_0 \longrightarrow \mathfrak{X}_{G_0}$ is
strictly singular.
\end{Proposition}
\noindent{\bf{Proof:}} Assume on the contrary. Then there exists a
block sequence $(x_n)_{n \in \N}$ with $\|x_n\|_{W_0}=1,n \in \N$
such that the operator \[ id|_{<x_n,n \in \N>}:(<x_n,n \in
\N>,\|\cdot\|_{W_0}) \to (<x_n,n \in \N>,\|\cdot\|_{G_0}) \] is an
isomorphism. Hence there exist $m>0$ and $M>0$ such that
\begin{equation}\label{L.5}
m\|x\|_{W_0} \le \|x\|_{G_0}\le M\|x\|_{W_0}, \mbox{for every~}
x \in <x_n,n \in \N>.
\end{equation}
Let $j\in \N$ with $\frac{1}{m_{2j+1}}<\frac{m}{90}$.
From $(\ref{L.5})$ we get that
$m\le\|x_n\|_{G_0}\le \|x_n\|_{W_0} \le 1$ for all $n \in \N$.
Therefore from Proposition $\ref{L.9}$ there exists a
$(15,2j+1,\frac{1}{2})$ dependent sequence
$(w_i,f_i)_{i=1}^{n_{2j+1}}$ in $\mathfrak{X}_0$
with $\ran(f_i)=\ran(w_i)$ and
$w_i\in <x_n,n \in \N>$ for all $i=1,\ldots,n_{2j+1}$.\\ We
set \[ w=\frac{m_{2j+1}}{n_{2j+1}}\sum_{i=1}^{n_{2j+1}}w_i
\mbox{~and~}
f=\frac{1}{m_{2j+1}}\sum\limits_{i=1}^{n_{2j+1}} f_i. \]
It is obvious that $f \in W_0$ and from the definition of
exact pairs in $\mathfrak{X}_0$
we get that $\|w\|_{W_0}\ge f(w)>\frac{1}{2}$.\\
On the other hand Proposition $\ref{L.4}$ yields
that $\|w\|_{G_0} \le\frac{45}{m_{2j+1}}.$ \\
Therefore from inequality $(\ref{L.5})$ it follows that $m\le
\frac{90}{m_{2j+1}}$, a contradiction.

\begin{Lemma}\label{L.3}
Let $\epsilon>0,\ (\epsilon_k)_{k\in \N}$ be a sequence of
positive numbers with $\sum\limits_{k=1}^{\infty}
\epsilon_k<\epsilon$, $(j_k)_{k\in\N}$ be a strictly increasing
sequence of even numbers and $(x_k)_{k \in \N}$ be a block
sequence such that
\begin{enumerate}
\item [i.] $\|x_k\|_{W_0}=1$ for all $k\in \N$
\item [ii.] $\frac{|\supp x_k|}{m_{j_k}}<\epsilon_k$,
for all $k\in\N$ and
\item [iii.] $\|x_{k+1}\|_{G_0} \le
\frac{\epsilon_k}{n_{j_k}}$ for
all $k\in\N$.
\end{enumerate}
Then $(x_k)_{k \in \N}$ is $(1+\epsilon)-$ equivalent to
the usual basis of $c_0(\N)$.
\end{Lemma}
\noindent{\bf{Proof:}}
Let $n \in \N$ and $a_1,\ldots,a_n \in \R$ with
$max\{|a_i|,i=1,\ldots,n\}=1$ for $i=1,\ldots,n$.
We assume without loss of generality that $n\ge 2$.
\\ Since the block sequence  $(x_k)_{k \in \N}$
is bimonotone and normalized it follows that
$\|\sum\limits_{i=1}^{n} a_ix_i\|_{W_0} \ge 1.$ We will
prove that
$\|\sum\limits_{i=1}^{n} a_ix_i\|_{W_0}\le 1+\epsilon$.
\\Let $f\in W_0$. We distinguish the following cases.\\
\textbf{Case 1.} Let $f \in G_0$. Then
\[|f(\sum\limits_{i=1}^{n} a_i x_i)| \le
\sum\limits_{i=1}^{n} |a_i|\cdot |f(x_i)| \le |f(x_1)|+
\sum\limits_{i=2}^{n} |f(x_i)| \le 1+
\sum\limits_{i=2}^{n} \epsilon_{i-1} \le 1+\epsilon.\]
\textbf{Case 2.} Let $f \in (W_0 \backslash G_0)$.\\
We assume without loss of generality that $f$ is of the form
$f=\frac{1}{m_{2j+1}}\sum\limits_{i=1}^{n_{2j+1}} f_i$, where
$(f_i)_{i=1}^{n_{2j+1}}$ is a $\sigma-n_{2j+1}$ special sequence.
Let $j_{k-1}<2j+1<j_k$ where $2\le k<n$. Then
\begin{eqnarray*}
|f(\sum\limits_{i=1}^{n} a_i x_i)| & & \le
|f(\sum\limits_{i=1}^{k-1} a_i x_i)|
+|f(a_k x_k)|+ |f(\sum\limits_{i=k+1}^{n} a_i x_i)|\le
\sum\limits_{i=1}^{k-1} |f(x_i)|+1+
\sum\limits_{i=k+1}^{n} |f(x_i)| \\ & &
\le 1+\sum\limits_{i=1}^{k-1} \frac{|supp(x_i)|}{m_{j_i}}+
\sum\limits_{i=k+1}^{n} n_{j_k-1}\cdot \|x_i\|_{G_0}=1+
\sum\limits_{i=1}^{k-1} \epsilon_i+\sum\limits_{i=k+1}^{n}
\frac{n_{j_k-1}\cdot \epsilon_i}{n_{j_{i-1}}} \\ & &
\le 1+\epsilon.
\end{eqnarray*}

\begin{Lemma}\label{L.52}
Let $\epsilon>0$ and $(x_k)_{k \in \N}$ be a block
sequence such that
\begin{enumerate}
\item [i.] $\lim\limits_{k \longrightarrow \infty}
\|x_k\|_{G_0}=0$ and
\item [ii.]  $\|x_k\|_{W_0}=1$ for all $k\in \N$.
\end{enumerate}
Then there exists a subsequence which is $(1+\epsilon)-$
equivalent to the usual basis of $c_0$.
\end{Lemma}
\noindent{\bf{Proof:}}
Let $(\epsilon_k)_{k\in \N}$ be a sequence of
positive numbers such that
$\sum\limits_{k=1}^{\infty} \epsilon_k<\epsilon$.
We inductively construct a subsequence
$(z_k)_{k \in \N}$ of $(x_k)_{k \in \N}$ and a strictly
increasing sequence $(j_k)_{k \in \N}$ of even numbers such
that
\begin{enumerate}
\item [i.] $\frac{|\supp z_k|}{m_{j_k}}<\epsilon_k$,
for all $k\in\N$ and
\item [ii.] $\|z_{k+1}\|_{G_0} \le
\frac{\epsilon_k}{n_{j_k}}$ for
all $k\in\N$.
\end{enumerate}
Therefore from Lemma $\ref{L.3}$, it follows that $(z_k)_{k \in
\N}$ is $(1+\epsilon)-$ equivalent to the usual basis of $c_0$.

\begin{Corollary}\label{L.102}
Let $Y$ be a subspace of $\mathfrak{X}_0$ and $\epsilon>0$. Then
there exists a subspace of $Y$ which is $(1+\epsilon)-$ isomorphic
to $c_0$.
\end{Corollary}
\noindent{\bf{Proof:}} We assume that $Y$ is a block subspace of
$\mathfrak{X}_0$. Then from the well known gliding
hump argument we get the result.\\
Let $(\epsilon_k)_{k\in \N}$ be a sequence of
positive numbers such that
$\sum\limits_{k=1}^{\infty} \epsilon_k<\epsilon$.\\
Since the identity operator $id:\mathfrak{X}_0 \longrightarrow
\mathfrak{X}_{G_0}$ is strictly singular, there exists a block
sequence $(y_k)_{k \in \N}$ in $Y$ such that
\begin{enumerate}
\item [i.] $\lim\limits_{k \longrightarrow \infty}
\|y_k\|_{G_0}=0$ and
\item [ii.]  $\|y_k\|_{W_0}=1$ for all $k\in \N$.
\end{enumerate}
Hence from Lemma $\ref{L.52}$ we obtain a subsequence which is
$(1+\epsilon)-$ equivalent to the usual basis of $c_0$.

\begin{Corollary}\label{L.54}
Let $Y$ be a subspace of $\mathfrak{X}_0$. Then $Y$ contains a
complemented copy of $c_0$.
\end{Corollary}
\noindent{\bf{Proof:}} Since $Y$ is separable and contains an
isomorphic copy of $c_0(\N)$, then from Sobczyk's \cite{SO}
theorem it follows that $c_0$ is complemented in $Y$.

\begin{Proposition}\label{L.50}
The basis $(e_n)_{n\in \N}$ of $\mathfrak{X}_0$ is shrinking,
hence $\mathfrak{X}_0^*$ is separable.
\end{Proposition}
\noindent{\bf{Proof:}} Suppose not. Then there exist
$\epsilon_0>0, x^* \in \mathfrak{X}_0^*$ with $\|x^*\|=1$ and a
block sequence $(x_n)_{n \in \N}$ with $\|x_n\|_{W_0} \le 1,n \in
\N$ such that $\epsilon_0<x^*(x_n)$ for all $n\in \N$. We
distinguish the following cases.
\vspace{1.5mm}\\
\textbf{Case 1.} There exists $\delta>0$ such that
$\delta<\|x_n\|_{G_0}$ for all $n \in \N$.
\vspace{2mm}\\
Proposition $\ref{L.24}$ yields that there exists
$n_0 \in \N$ such that for every $k \in \N$ with
$k>n_0$ and every $L \in [\N]$ there exists
$Q \in [L]$ such that
\\ for every $n_1<\ldots<n_k$ in $Q$ we have that
\[\|\frac{x_{n_1}+\ldots+x_{n_k}}{k}\|_{W_0}=
\|\frac{x_{n_1}+\ldots+x_{n_k}}{k}\|_{G_0}.\]
Hence there exists $M \in [\N]$ and a block sequence
$(z_k)_{k \in \N}$ of $(x_k)_{k \in \N}$ such that
\begin{enumerate}
\item $z_k=\frac{1}{|F_k|} \sum\limits_{i \in F_k} x_i$,
where $F_k<F_{k+1}$ finite subsets of $M$ and
$(|F_k|)_{k \in \N}$ strictly increasing sequence with
$|F_1|>n_0$
\item $\epsilon_0<x^*(z_k)<\|z_k\|_{W_0}=\|z_k\|_{G_0}
\le 1$ for all $k \in \N$.
\end{enumerate}
Let $j \in \N$ with $n_{2j}>k_0$ and $\epsilon>0$ with
$\epsilon<\frac{2}{m_{2j}^2}$.
We may assume without loss of generality that
$(z_k)_{k \in \N}$ is $(3,\epsilon)$ RIS in
$\mathfrak{X}_0$ such that the first term of the
associated sequence is bigger than $2j$.\\
From (2) it follows that there exist
$k_1<\ldots<k_{n_{2j}}$ in $M$ such that
\[\|\frac{z_{k_1}+\ldots+z_{k_{n_{2j}}}}{n_{2j}}\|_{W_0}=
\|\frac{z_{k_1}+\ldots+z_{k_{n_{2j}}}}{n_{2j}}\|_{G_0}.\]
Thus from proposition $\ref{L.8}$ we have that
\[\|\frac{z_{k_1}+\ldots+z_{k_{n_{2j}}}}{n_{2j}}\|_{W_0}=
\|\frac{z_{k_1}+\ldots+z_{k_{n_{2j}}}}{n_{2j}}\|_{G_0} \le
\frac{9}{m_{2j}}.\] Since the action of $x^*$ in every convex
combination of $(z_k)_{k \in \N}$ is bigger that $\epsilon_0$,
then for a sufficiently large $j \in \N$ we derive to
contradiction.
\vspace{2mm}\\
\textbf{Case 2.} There exists a subsequence
$(x_{p_k})_{k \in \N}$ such that
$\lim\limits_{k \longrightarrow \infty} \|x_{p_k}\|_{G_0}=0$.\\
\vspace{2mm} Then from Lemma $\ref{L.3}$ we assume without loss of
generality that $(x_{p_k})_{k \in \N}$ is equivalent to the usual
basis of $c_0$. Hence there exists $M>0$ such that
\[k\epsilon_0 \le x^*(\sum\limits_{i=1}^{k} x_{p_i}) \le
\|\sum\limits_{i=1}^{k} x_{p_i}\|_{W_0} \le M
\mbox{~for every~} k \in \N\] a contradiction.

\vspace{4mm} Next we shall show a structural property for the
subspaces of $\mathfrak{X}_0$ on which the identity operator
$id:\mathfrak{X}_0 \longrightarrow \mathfrak{X}_{G_0}$ is compact.
Our approach uses a beautiful and deep result due to Kalton
\cite{K}.

\begin{Notation}
We set $[\N]^{<\omega}=\bigcup_{k \in \N} [N]^k$, where $[\N]^k$
denotes the set of all finite subsets of $\N$ of cardinality $k$.
We define in $[\N]^{<\omega}$ the
following partial order: \\
if $\{n_1<\ldots<n_k\},\{m_1<\ldots<m_n\} \in [\N]^{<\omega}$
then
\[ \{n_1<\ldots<n_k\} \sqsubseteq \{m_1<\ldots<m_n\}
\mbox{~if and only if~} k \le n \mbox{~and~} n_i=m_i \mbox{~for
all~} i=1,\ldots,k.\] This partial order is called the initial
segment partial order and the couple
$([\N]^{<\omega},\sqsubseteq)$ is a tree. A branch of this tree is
identified by an infinite subset $\{q_n,n \in \N\}$ of natural
numbers, where $(q_n)_{n \in \N}$ is a strictly increasing
sequence. Moreover if $b=\{k_1<\ldots<k_n<\ldots\}$ is a branch of
this tree, then the node $\{k_1<\ldots<k_n\}$ is symbolized as
$b|n$.
\\ Let $k \in \N$ and $s \in [\N]^{<\omega}$ with
$s=\{m_1<\ldots<m_k\}$. If $m \in \N$ with $m_k<m$, then the
finite subset $\{m_1<\ldots<m_k<m\}$, is symbolized as $(s^ \frown
m)$.
\end{Notation}

\begin{Definition}
Let $(X,\|\cdot\|)$ be an infinite dimensional Banach space.
A family of vectors $(x_s)_{s \in [\N]^{<\omega}}$ in $X$
is called normalized weakly null tree family in $X$ if:
\begin{enumerate}
\item [i.] $\|x_s\|=1$ for all $s \in [\N]^{<\omega}$ and
\item [ii.] for each node $s \in [\N]^{<\omega}$ the
sequence $(x_{(s ^ \frown n)})_{n \in \N}$ is weakly null.
\end{enumerate}
\end{Definition}

\begin{Definition}\label{L.120}
Let $(X,\|\cdot\|)$ be an infinite dimensional Banach space.
We say that $X$ has the $c_0$ tree property if
there exists $K>0$ such that for every normalized weakly null
tree family $(x_s)_{s \in [\N]^{<\omega}}$ in $X$ there exist
a branch $b$ of the tree $([\N]^{<\omega},\sqsubseteq)$ such
that the sequence $(x_{b|n})_{n \in \N}$ is $K-$
equivalent to the usual basis of $c_0(\N)$.
\end{Definition}

\vspace{1mm} We pass to state Kalton's Theorem (\cite{K}, Thm.
3.2).

\begin{Theorem}\label{L.36}
Let $X$ be a separable Banach space not containing $\ell^1(\N)$.
If $X$ has the $c_0$ tree property, then $X$ is embedded into
$c_0$.
\end{Theorem}

\begin{Remark}
Kalton's theorem provides an efficient characterization of the
subspaces of $c_0$. Similar results for $\ell^p$ spaces have been
proved by Odell and Schlumprecht \cite{OS2} and for subspaces of
reflexive spaces with an unconditional basis by W. B. Johnson and
B. Zheng \cite{JZ}.
\end{Remark}

\begin{Proposition}\label{L.11}
Let $Y$ be a subspace of $\mathfrak{X}_0$ such that the identity
operator $id|_{Y}:(Y,\|\cdot\|_{W_0}) \longrightarrow
\mathfrak{X}_{G_0}$ is compact. Then $Y$ has the $c_0$ tree
property, hence $Y$ is embedded isomorphically into $c_0$.
\end{Proposition}
\noindent{\bf{Proof:}}
Let $\epsilon>0$ and $(x_a)_{a \in [\N]^{<\omega}}$ be
a normalized weakly null tree family in
$(Y,\|\cdot\|_{W_0})$.
We will construct by induction two strictly increasing
sequences $(l_k)_{k\in \N}$ and
$(j_k)_{k \in \N}$ of natural numbers and a block
sequence $(z_k)_{k\in \N}$ such that:
\begin{enumerate}
\item [(i)] $\frac{1}{2}\le\|z_k\|_{W_0}\le 2$ for all $k \in \N$
\item [(ii)] $\frac{|\supp z_k|}{m_{j_k}}<\epsilon_k$, for all
$k\in\N$ \item [(iii)] $\|z_k\|_{G_0} \le
\frac{\epsilon_{k-1}}{n_{j_{k-1}}}$ for all $k\in\N$ with $k\ge 2$
\item [(iv)] setting $s_k=\{l_1<\ldots<l_k\}$ for all $k \in \N$,
then $\|x_{s_k}-z_k\|_{W_0}<\epsilon_k$ for all $k \in \N$,
\end{enumerate}
where $(\epsilon_k)_{k \in \N}$ is a sequence of
positive numbers with $\epsilon_k<\frac{1}{2}$ for all
$k \in \N$
and $\sum\limits_{k=1}^{\infty} \epsilon_k<\frac{\epsilon}{8}$.
\\ Indeed, let $l_1 \in \N$. Since $x_{l_1} \in \mathfrak{X}_0$
there exists a block vector $z_1$ such that
$\|x_{l_1}-z_1\|_{W_0}<\epsilon_1$. We choose $j_1 \in \N$ such
that $\frac{|\supp z_1|}{m_{j_1}}<\epsilon_1$.
It is obvious that condition (i) is satisfied.\\
We assume that for a fixed $k \in \N$ with $k\ge 2$ we have
constructed a finite increasing sequence $(l_i)_{i=1}^{k}$
of natural numbers,
a finite sequence $(j_i)_{i=1}^{k}$ of increasing natural
numbers
and a finite block sequence $(z_i)_{i=1}^{k}$ such that are
satisfied conditions (i)-(iv) until the fixed $k \in \N$.\\
The sequence $(x_{(s_k^\frown n)})_{n \in \N}$ is weakly null and
normalized, therefore there exists a subsequence $(x_{(s_k ^\frown
q_n)})_{n \in \N}$ of $(x_{(s_k ^\frown n)})_{n \in \N}$ and a
block sequence $(z_n^k)_{n \in \N}$ such that $\|x_{(s_k ^\frown
q_n)}-z_n^k\|_{W_0}<r_n$ for all $n \in \N$, where
$(r_n)_{n \in \N}$ is a null sequence of positive numbers.\\
Since the sequence $(x_{(s_k ^\frown q_n)})_{n \in \N}$ is
weakly null and the identity operator
$id:(Y,\|\cdot\|_{W_0}) \longrightarrow \mathfrak{X}_{G_0}$ is
compact, it follows that $\lim\limits_{n} \|z_n^k\|_{G_0}=0$.\\
Therefore there exists $n_k \in \N$ with $q_{n_k}>l_k$ such that
$\|z_{n_k}^k\|_{G_0}<\frac{\epsilon_k}{n_{j_k}},z_k<z_{n_k}^k$
and $\|x_{(s_k ^\frown q_{n_k})}-z_{n_k}^k\|_{W_0}<\epsilon_{k+1}.$
\\ We set $z_{k+1}=z_{n_k}^k$ and $l_{k+1}=q_{n_k}$.
We choose $j_{k+1} \in \N$ with
$j_k<j_{k+1}$ such that
$\frac{|\supp z_{k+1}|}{m_{j_{k+1}}}<\epsilon_{k+1}$.\\
Thus we have constructed a block sequence $(z_k)_{k \in \N}$ and a
branch $\{l_1<\ldots<l_k<\ldots\}$ of the tree $[\N]^{<\omega}$,
such that the sequence $(x_{s_k})_{k \in \N}$ is 4-equivalent to
$(z_k)_{k \in \N}$. Conditions (i), (ii), (iii) yield that the
block sequence $(z_k)_{k \in \N}$ is $4+\frac{\epsilon}{4}$
equivalent to the usual basis $(e_k)_{k \in \N}$ of $c_0$(Lemma
$\ref{L.3}$) and therefore $(x_{s_k})_{k\in \N}$ is $16+\epsilon$
equivalent to
the usual basis $(e_k)_{k \in \N}$ of $c_0$.\\
Also, since $Y$ is separable, does not contain an isomorphic copy
of $\ell^1(\N)$ and has the $c_0$ tree property, it follows from
Thm. $\ref{L.36}$ that $Y$ is embedded isomorphically into $c_0$.

\begin{Proposition}\label{L.49}
Let $Y,Z$ be subspaces of $\mathfrak{X}_0$ such that the identity
operators $id|_{Y}:(Y,\|\cdot\|_{W_0}) \longrightarrow
\mathfrak{X}_{G_0}$, $id|_{Z}:(Z,\|\cdot\|_{W_0}) \longrightarrow
\mathfrak{X}_{G_0}$ are not compact. Then $d(S_Y,S_Z)=0$.
\end{Proposition}
\noindent{\bf{Proof:}}
We shall show that for every $\epsilon>0$ there exist
$y \in Y,z \in Z$ with
$\|z-y\|_{W_0}<\epsilon \cdot \|z+y\|_{W_0}$.\\
Then as is known this yields the result.\\
From the fact that $id|_{Y}:(Y,\|\cdot\|_{W_0}) \longrightarrow
\mathfrak{X}_{G_0}$ is not compact and $Y$ does not contain an
isomorphic copy of $\ell^{1}(\N)$, it follows that there exist a
$\|\cdot\|_{W_0}$ normalized sequence $(y_n)_{n \in \N}$ in $Y$
and a block sequence $(x_n)_{n \in \N}$ in $\mathfrak{X}_0$ such
that these sequences are equivalent and the identity operator
$id|_{(\overline{span}\{x_n:n \in \N\},\|\cdot\|_{W_0})}$ is not
compact. The same holds for the subspace $Z$. Hence using the well
known gliding hump argument, we may assume
that $Y,Z$ are block subspaces.\\
Since the identity operators are not compact, it follows that
there exist block sequences $(y_n)_{n \in \N}$ in $Y$, $(z_n)_{n
\in \N}$ in $Z$ and $\epsilon_0>0$ such that
\[\epsilon_0<\|y_n\|_{G_0} \le \|y_n\|_{W_0}\le 1
\mbox{~for all~} n \in \N \]
and
\[\epsilon_0<\|z_n\|_{G_0} \le \|z_n\|_{W_0}\le 1
\mbox{~for all~} n \in \N.\] Therefore, Proposition $\ref{L.20}$
yields that for every $j \in \N$ there exists a pair $(w,f)$ with
$w \in <y_n:n \in \N>$ or $w \in <z_n:n \in \N>$ which is
$(15,2j,\frac{1}{2})$ exact pair in
$\mathfrak{X}_0$.\\
Hence if $j_0 \in \N$, then there exists a finite sequence of
pairs $(w_k,f_k)_{k=1}^{n_{2j_0+1}}$ which is
$(15,2j_0+1,\frac{1}{2})$ dependent sequence in $\mathfrak{X}_0$
with $\ran(f_k)=\ran(w_k)$ for all $k=1,\ldots,n_{2j_0+1}$ such
that for every $k \in \{1,\ldots,n_{2j_0+1}\}$ we have that $w_k
\in <y_n:n \in \N>$ if $k$ is odd and $w_k \in <z_n:n \in \N>$ if
$k$ is even
(Prop. $\ref{L.9}$).\\
From Proposition $\ref{L.2}$ we get that
\begin{equation}\label{L.29}
\|\frac{1}{n_{2j_0+1}} \sum\limits_{i=1}^{n_{2j_0+1}}
(-1)^{i}w_i\|_{W_0} \le \frac{120}{m_{2j_0+1}^2}.
\end{equation}
Also we have that
\begin{equation}\label{L.30}
\|\frac{1}{n_{2j_0+1}}
\sum\limits_{i=1}^{n_{2j_0+1}} w_i\|_{W_0} \ge
(\frac{1}{m_{2j_0+1}}\sum\limits_{i=1}^{n_{2j_0+1}} f_i)
(\frac{1}{n_{2j_0+1}}
\sum\limits_{i=1}^{n_{2j_0+1}} w_i) \ge
\frac{1}{2m_{2j_0+1}}.
\end{equation}
We set \[A=\{i\in \{1,\ldots,n_{2j_0+1}\}:i \mbox{~even} \}
\mbox{~and~}
B=\{i\in \{1,\ldots,n_{2j_0+1}\}:i \mbox{~odd} \}.\]
Moreover setting $y=\frac{1}{n_{2j_0+1}}\sum\limits_{i \in B} w_i$
and
$z=\frac{1}{n_{2j_0+1}}\sum\limits_{i \in A} w_i$,
from $(\ref{L.29})$ and $(\ref{L.30})$,
we get that
\[\|z-y\|_{W_0}<\frac{240}{m_{2j_0+1}}\cdot
\|z+y\|_{W_0}\]
which implies that
$\|z-y\|_{W_0}<\epsilon\cdot \|z+y\|_{W_0}$ for fixed
$\epsilon>0$ and a sufficiently large $j_0 \in \N$.

\begin{Corollary}\label{L.38}
Let $Y,Z$ be subspaces of $\mathfrak{X}_0$, such that the direct
sum $Y\bigoplus Z$ is a closed subspace of $\mathfrak{X}_0$. Then
at least one of the identity operators
$id|_{Y}:(Y,\|\cdot\|_{W_0}) \longrightarrow \mathfrak{X}_{G_0}$,
$id|_{Z}:(Z,\|\cdot\|_{W_0}) \longrightarrow \mathfrak{X}_{G_0}$
is compact.\\ Hence one at least of the subspaces is embedded
isomorphically into $c_0$.
\end{Corollary}
\noindent{\bf{Proof:}} Assume on the contrary. Then from
Proposition $\ref{L.49}$ we get that $d(S_Y,S_Z)=0$, a
contradiction. Moreover from Proposition $\ref{L.11}$ we obtain
that one at least of the subspaces is embedded isomorphically into
$c_0$.

\begin{Remark}\label{L.44}
Let $Y,Z$ be subspaces of $\mathfrak{X}_0$ such that
$\mathfrak{X}_i=Y \bigoplus Z$. Then from the fact that the
identity operator $id:\mathfrak{X}_0 \longrightarrow
\mathfrak{X}_{G_0}$ is not compact, it follows that one at least
of the identity operators $id|_{Y}:(Y,\|\cdot\|_{W_0})
\longrightarrow \mathfrak{X}_{G_0}$ or
$id|_{Z}:(Z,\|\cdot\|_{W_0}) \longrightarrow \mathfrak{X}_{G_0}$
is not compact. Hence Corollary $\ref{L.38}$ yields that $id|_{Y}$
is compact and $id|_{Z}$ is not compact or vice versa.
\end{Remark}

\begin{Corollary}\label{L.45}
Let $Y$ be a subspace of $\mathfrak{X}_0$ such that the identity
operator $id|_{Y}:(Y,\|\cdot\|_{W_0}) \longrightarrow
\mathfrak{X}_{G_0}$ is not compact. Then $(Y,\|\cdot\|_{W_0})$ is
not embedded isomorphically into $c_0$.
\end{Corollary}
\noindent{\bf{Proof:}} Suppose not. Then $(Y,\|\cdot\|_{W_0})$
embeds isomorphically into $c_0(\N)$. Since the operator $id|_{Y}$
is not compact and $Y$ does not contain an isomorphic copy of
$\ell^{1}(\N)$, it follows that there exist a sequence $(y_n)_{n
\in \N}$ in $Y$ and a block sequence $(x_n)_{n \in \N}$ in
$\mathfrak{X}_0$ such that these sequences are equivalent and the
identity operator $id|_{(\overline{span}\{x_n:n \in \N\},
\|\cdot\|_{W_0})}$ is not compact. Hence there exist a block
sequence $(w_n)_{n \in \N}$ of $(x_n)_{n \in \N}$ and $\epsilon>0$
such that
\[ \epsilon<\|w_n\|_{G_0}\le \|w_n\|_{W_0}\le 1 \mbox{~for all~}
n \in \N.\] Since the basis of $\mathfrak{X}_0$ is shrinking, we
obtain that $(w_n)_{n \in \N}$ is weakly null in $\mathfrak{X}_0$.
From the fact that $(\overline{span}\{x_n:n \in
\N\},\|\cdot\|_{W_0})$ is embedded isomorphically into $c_0$, we
get (without loss of generality) that $(w_n)_{n \in \N}$ is
equivalent to the usual basis of $c_0(\N)$. Hence there exist
$m,M>0$ such that
\begin{equation}\label{L.43}
m\cdot max\{|a_i|,i=1,\ldots,n\} \le \|\sum\limits_{i=1}^{n}
a_iw_i\|_{W_0}\le M\cdot max\{|a_i|,i=1,\ldots,n\}
\end{equation}
for all $n \in \N$ and $a_1,\ldots,a_n$ real numbers.\\
For every $n \in \N$ there exists $f_n \in G_0$ such that
$ran(f_n) \subset ran(w_n)$ and $f_n(w_n)>\epsilon$.\\
Thus setting
\[g_j=\frac{1}{m_{2j}} \sum\limits_{i=1}^{n_{2j}} f_i,
j \in \N\]
we get that $g_j \in G_0$ and
$\epsilon \le \|\frac{m_{2j}}{n_{2j}}
\sum\limits_{i=1}^{n_{2j}} w_i\|_{G_0}$ for all $j \in \N$.
From $(\ref{L.43})$ we have that \\
$\|\frac{m_{2j}}{n_{2j}} \sum\limits_{i=1}^{n_{2j}} w_i\|_{G_0}
\le M\cdot \frac{m_{2j}}{n_{2j}}$, for all $j \in \N$, a
contradiction for a sufficiently large $j \in \N$.
\vspace{2.5mm}\\
We pass to the definition of type I and type II
complemented subspaces of $\mathfrak{X}_0$.

\begin{Definition}
Let $Y$ be a complemented subspace of $\mathfrak{X}_0$. We say
that \item [i.] $Y$ is of type I if contains an isomorph of
$\mathfrak{X}_0$. \item [ii.] $Y$ is of type II if it is
isomorphic to a subspace of $c_0$.
\end{Definition}

\begin{Theorem}\label{L.41}
Let $Y,Z$ be subspaces of $\mathfrak{X}_0$ such that
$\mathfrak{X}_0=Y\bigoplus Z$. Then $Y$ is of type I and $Z$ is of
type II or vice versa. If especially $Y\cong c_0$, then $Z\cong
\mathfrak{X}_0$.
\end{Theorem}
\noindent{\bf{Proof:}} From Corollary $\ref{L.38}$, Remark
$\ref{L.44}$ and Corollary $\ref{L.45}$ it follows that $Y$ is
embedded isomorphically into $c_0$ and $Z$ is not embedded
isomorphically into $c_0(\N)$ or vice versa. Assume without loss
of generality the first case. Then there exists a subspace $Z_1$
of $Z$ such that $Z\cong c_0 \bigoplus Z_1$. Hence
\[\mathfrak{X}_0=Y\bigoplus Z \cong Y \bigoplus c_0
\bigoplus Z_1\] and since $Y$ is embedded isomorphically into
$c_0$ and $c_0\bigoplus c_0 \cong c_0$, we get the conclusion.\\
In the special case where $Y\cong c_0$ we have that
\[\mathfrak{X}_0=Y\bigoplus Z \cong c_0 \bigoplus
c_0 \bigoplus Z_1 \cong c_0 \bigoplus Z_1 \cong Z.\]

\vspace{3mm} P. Koszmider \cite{KO}, under CH, has constructed a
nonseparable $C(K)$ space satisfying the property that whenever
$Y,Z$ be subspaces of $C(K)$ such that $C(K)=Y \bigoplus Z$, then
either $Y \cong c_0$ and $Z \cong C(K)$ or vice versa. In the same
paper he asked whether a separable Banach space could occur
sharing similar properties. The answer of this problem is
affirmative and was given by S. A. Argyros and Th. Raikoftsalis in
\cite{AR}. In this paper it is introduced a new class of primary
Banach spaces called quasi-prime. An infinite dimensional Banach
space $X$ is said to be primary, if $Y,Z$ be closed subspaces of
$X$ such that $X=Y\bigoplus Z$, then $Y\cong X$ or $Z \cong X$.
The quasi-prime Banach spaces are spaces which satisfy a property
like the above $C(K)$. In the present paper, the Banach space
$\mathfrak{X}_0$, resembles the quasi-prime Banach spaces as it
seems in Thm. $\ref{L.41}$.

\begin{Corollary}
The basis $(e_n)_{n \in \N}$ of $\mathfrak{X}_0$, is a normalized,
weakly null sequence without unconditional subsequence.
\end{Corollary}
\noindent{\bf{Proof:}} Assume the contrary. Then there exists an
$L \in [\N]$ such that $(e_n)_{n \in L}$ is unconditional. Let
$L_1,L_2$ in $[L]$ with $L_1 \cap L_2=\emptyset$ and $L=L_1\cup
L_2.$ Then
\[ \overline{<e_n,n \in L>}=\overline{<e_n,n \in L_1>}
\bigoplus \overline{<e_n,n \in L_2>}.\] From Corollary
$\ref{L.38}$ we obtain that one at least of the identity operators
$id|_{\overline{<e_n,n \in L_1>}},id|_{\overline{<e_n,n \in
L_2>}}$ is compact. Let $id|_{\overline{<e_n,n \in L_1>}}$ is
compact. Since the basis of $\mathfrak{X}_0$ is shrinking, it
follows that $(e_n)_{n \in L_1}$ is weakly null in
$\mathfrak{X}_0$ and from the compactness of the identity operator
we get that $\lim\limits_{n \in L_1} \|e_n\|_{G_0}=0$, a
contradiction.

\begin{Remark}
The space $\mathfrak{X}_0$ is not embedded into
a space with an unconditional basis.
\end{Remark}

\section{\textbf{The space $\mathfrak{L}(\mathfrak{X}_0)$}}

In this section we study the structure of the operators of the
space $\mathfrak{X}_0$. Since the space $\mathfrak{X}_0$ is $c_0$
saturated it admits many projections. The aim is to show Theorem
$\ref{L.22}$, which asserts that beyond the identity the non
strictly singular operators are isomorph only on subspaces of
$\mathfrak{X}_0$ which are embedded into $c_0$. \vspace{2mm}

\begin{Notation}
Let $x \in \mathfrak{X}_0$ and $Y$ be a subspace of
$\mathfrak{X}_0$. We note by $d_{G_0}(x,Y)$ the following number
\[d_{G_0}(x,Y)=inf\{\|x-y\|_{G_0},y \in Y\}.\]
\end{Notation}

\begin{Lemma}\label{L.27}
Let $(e_n)_{n \in \N}$ be the basis of $\mathfrak{X}_{G_0}$ and
$\mathfrak{X}_0$ and
$T:\mathfrak{X}_0 \longrightarrow \mathfrak{X}_0$ be a bounded
linear operator.\\ Then
$\lim\limits_{n \longrightarrow \infty}
d_{G_0}(T(e_n),\R e_n)=0.$
\end{Lemma}
\noindent{\bf{Proof:}} Assume on the contrary. Then there exist
$\delta>0$ and an infinite subset $L$ of $\N$ such that
\begin{equation}\label{L.47}
0<\delta\le d_{G_0}(T(e_n),\R e_n)=
inf\{\|T(e_n)-x\|_{G_0},x \in \R e_n\}
\mbox{~for all~} n \in L.
\end{equation}
Since $\delta \le \|T(e_n)\|_{G_0}$ for all $n \in L$ and
$(T(e_n))_{n \in \N}$ is weakly null in $\mathfrak{X}_0$
we may assume without loss of generality that there exist
a block sequence $(x_n)_{n \in L}$ such that
\begin{equation}\label{L.48}
\|T(e_n)-x_n\|_{W_0}<\epsilon_n \mbox{~for all~} n \in L
\end{equation}
where $(\epsilon_n)_{n \in L}$ is a sequence of positive
numbers with $\epsilon_n<\frac{\delta}{2}$ for all
$n \in L$ and
$\sum\limits_{n=1}^{\infty} \epsilon_n \le 1$.
\\ For every $n \in L$ we split $x_n$ into tree vectors
as follows
\[x_n=x_n^{(1)}+x_n^{(2)}+x_n^{(3)} \mbox{~with~}
x_n^{(i)}=E_n^{(i)}x_n  \mbox{~for~} i=1,2,3\] where
$E_n^{(1)}=\{k \in \N:k<n\},E_n^{(3)}=\{k \in \N:k>n\}$
and $E_n^{(2)}=\{n\}$.\\
From inequalities $(\ref{L.47})$ and $(\ref{L.48})$
we conclude that for every $n \in L$ either
$\frac{\delta}{4}<\|x_n^{(1)}\|_{G_0}$ or
$\frac{\delta}{4}<\|x_n^{(3)}\|_{G_0}$.
Hence we may assume without loss of generality that
$\frac{\delta}{4}<\|x_n^{(1)}\|_{G_0}$ for all
$n \in L$ or $\frac{\delta}{4}<\|x_n^{(3)}\|_{G_0}$
for all $n \in L$.\\
Let $\frac{\delta}{4}<\|x_n^{(1)}\|_{G_0}$ for all $n \in L$. In
the other case the proof is similar as
follows.\\
Since $(x_n^{(1)})_{n \in L_0}$ is block, we get that for every $n
\in L_0$ there exists $g_n \in G_0$ with $ran(g_n)\subset
ran(x_n^{(1)})$ and $g_n(x_n^{(1)})>\frac{\delta}{4}$. Moreover we
observe that $g_n(x_n^{(1)})=g_n(x_n)$ for all
$n \in L_0$.\\
Let $j \in \N$ and $k_1<\ldots<k_{n_{2j}}$ natural numbers in
$L_0$ with $ran(x_{k_1})<k_1$ and $k_i<ran(x_{k_{i+1}})<k_{i+1}$
for $i=1,\ldots,n_{2j}-1$. Setting
\[y_j=\frac{m_{2j}}{n_{2j}}\sum\limits_{r=1}^{n_{2j}}
x_{k_r},z_j=\frac{m_{2j}}{n_{2j}}
\sum\limits_{r=1}^{n_{2j}} e_{k_r}
\mbox{~and~}
f_j=\frac{1}{m_{2j}}\sum\limits_{r=1}^{n_{2j}} g_{k_r}\]
we observe that
\begin{enumerate}
\item [i.] $\frac{\delta}{4}\le f_j(y_j)$ and
\item [ii.] $f_j(z_j)=0$.
\end{enumerate}
Also from $(\ref{L.48})$ we have that
\begin{equation}\label{L.35}
\|T(z_j)-y_j\|_{W_0}\le
\frac{m_{2j}}{n_{2j}}\sum\limits_{r=1}^{n_{2j}}\epsilon_r
\le \frac{1}{m_{2j}}.
\end{equation}
Let $j_0 \in \N$. Then we construct a finite sequence of pairs
$(z_j,f_j)_{j=1}^{n_{2j_0+1}}$ which is a $(15,2j_0+1,0)$
dependent sequence in $\mathfrak{X}_0$. Hence Proposition
$\ref{L.2}$ yields that
\begin{equation}\label{L.46}
\|\frac{m_{2j_0+1}}{n_{2j_0+1}}
\sum\limits_{j=1}^{n_{2j_0+1}} z_j\|_{W_0} \le
\frac{120}{m_{2j_0+1}}
\end{equation}
Therefore from inequality $(\ref{L.35})$ and $(\ref{L.46})$
we get that
\begin{eqnarray*}
\frac{\delta}{4}\le \|\frac{m_{2j_0+1}}{n_{2j_0+1}}
\sum\limits_{j=1}^{n_{2j_0+1}} y_j\|_{W_0} & & \le
\|\frac{m_{2j_0+1}}{n_{2j_0+1}} \sum\limits_{j=1}^{n_{2j_0+1}}
(y_j-T(z_j))\|_{W_0}+
\|\frac{m_{2j_0+1}}{n_{2j_0+1}}
\sum\limits_{j=1}^{n_{2j_0+1}}
T(z_j))\|_{W_0} \\ & & \le \frac{1}{m_{2j_0+1}}+
\|T\| \|\frac{m_{2j_0+1}}{n_{2j_0+1}}
\sum\limits_{j=1}^{n_{2j_0+1}} z_j\|_{W_0} \le
\frac{120\|T\|+1}{m_{2j_0+1}}
\end{eqnarray*}
which is a contradiction for a sufficiently large
$j_0 \in \N$.

\begin{Lemma}\label{L.28}
Let $T:\mathfrak{X}_0 \longrightarrow \mathfrak{X}_0$ be a
bounded linear operator and $L \in [\N]$ such that
\begin{enumerate}
\item [i.]  $\lim\limits_{n \in L} \|T(e_n)\|_{G_0}=0$ and
\item [ii.] there exists a subspace $Y$ of $\mathfrak{X}_0$
such that the restriction
$T_{|Y}:Y \longrightarrow T(Y)$ of $T$  is an
isomorphism.
\end{enumerate}
Then the identity operator
$id:(Y,\|\cdot\|_{W_0}) \longrightarrow \mathfrak{X}_{G_0}$
is compact.
\end{Lemma}
\noindent{\bf{Proof:}} Suppose not. Then as in Proposition
$\ref{L.49}$ we may assume without loss of generality that $Y$ is
a block subspace of $\mathfrak{X}_0$. Hence there exists a block
sequence $(y_n)_{n \in \N}$ in $Y$ and $\epsilon_0>0$ such that
\[\epsilon_0<\|y_n\|_{G_0}\le \|y_n\|_{W_0}\le 1
\mbox{~for all~} n \in \N.\] Since $T_{|Y}$ is an isomorphism,
there exists $m>0$ such that
\begin{equation}\label{L.32}
m \|x\|_{W_0} \le \|T(x)\|_{W_0} \le \|T\| \|x\|_{W_0}
\mbox{~for every~} x \in Y.
\end{equation}
We distinguish the following cases.\\
\vspace{1mm}\\
\textbf{1.} Let $\lim\limits_{n \in L}
\|T(e_n)\|_{W_0}=0$.\\
\vspace{1mm}\\
There exists $M \in [L]$ such that
\begin{equation}\label{L.56}
\sum\limits_{n \in M} \|T(e_n)\|_{W_0}<\frac{m}{2}.
\end{equation}
From the fact that the identity operators
$id|_{Y}:(Y,\|\cdot\|_{W_0}) \longrightarrow \mathfrak{X}_{G_0}$,
and $id|_{\overline{<e_n,n \in M>}}: (\overline{<e_n,n \in
M>},\|\cdot\|_{W_0}) \longrightarrow \mathfrak{X}_{G_0}$ are not
compact, Prop. $\ref{L.49}$ yields that there exist $y \in Y$ with
$\|y\|_{W_0}=1$ and $z \in <e_n,n \in M>$ with $\|z\|_{W_0}=1$
such that
$\|y-z\|_{W_0}<\frac{m}{2\|T\|}$.\\
Therefore
\begin{equation}\label{L.55}
\|T(y-z)\|_{W_0}\le \|T\| \cdot \|y-z\|_{W_0}
<\frac{m}{2}.
\end{equation}
Let $z=\sum\limits_{i \in F} a_ie_i$, where $F$ is a finite subset
of $M$. Then using $(\ref{L.56})$ we get that
\begin{equation}\label{L.57}
\|T(z)\|_{W_0} \le \sum\limits_{i \in F}
\|T(e_i)\|_{W_0}<\frac{m}{2}.
\end{equation}
Hence from the isomorphism and $(\ref{L.57})$
we have that
\begin{equation}\label{L.58}
\|T(y-z)\|_{W_0}\ge \|T(y)\|_{W_0}-\|T(z)\|_{W_0}
\ge m-\frac{m}{2}=\frac{m}{2}
\end{equation}
which contradicts to $(\ref{L.55})$.
\vspace{2mm}\\
\textbf{2.} Let $\limsup\limits_{n \in L }
\|T(e_n)\|_{W_0}>0$.\\
\vspace{2mm}\\
We set $L=\{l_1<\ldots<l_n<\ldots\}$. Since $(T(e_n))_{n \in \N}$
is weakly null in $\mathfrak{X}_0$, we may assume without loss of
generality that there exist a block sequence $(x_n)_{n \in L}$
which is equivalent to $(T(e_n))_{n \in L}$. Hence the block
sequence $(x_n)_{n \in L}$ is $\|\cdot\|_{W_0}-$ bounded,
$\|\cdot\|_{W_0}-$ seminormalized and
$\lim\limits_{n \in L} \|x_n\|_{G_0}=0$.\\
From Lemma $\ref{L.52}$ there exists a subsequence of $(x_n)_{n\in
L}$ which is equivalent to the usual basis of $c_0$. We may assume
without loss of generality that $(x_n)_{n \in L}$ is equivalent to
the usual basis of $c_0(\N)$. Moreover, since $(x_n)_{n \in L}$ is
equivalent to $(T(e_n))_{n \in L}$, it follows that there exist
$d_1,d_2>0$ such that
\begin{equation}\label{L.33}
d_1 \cdot max\{|a_i|:i=1,\dots,n\}\le
\|\sum\limits_{i=1}^{n} a_i T(e_{l_i})\|_{W_0}\le d_2
\cdot max\{|a_i|:i=1,\dots,n\}
\end{equation}
for every $n \in \N$ and $a_1,\ldots,a_n \in \R$. \\
We choose $j_0 \in \N$ such that
$\frac{1}{m_{2j_0+1}}<\min\{\frac{m}{8d_2},
\frac{m}{8 \cdot 120\|T\|}\}$.\\
Hence from Prop. $\ref{L.9}$ there exists a
$(15,2j_0+1,\frac{1}{2})$ dependent sequence
$(w_k,f_k)_{k=1}^{n_{2j_0+1}}$ in $\mathfrak{X}_0$ with
$\ran(f_k)=\ran(w_k)$ for all $k=1,\ldots,n_{2j_0+1}$ such that
for every $k \in \{1,\ldots,n_{2j_0+1}\}$ we have that $w_k \in
<y_n:n \in \N>$ if $k$ is odd and
$w_k \in <e_n:n \in \N>$ if $k$ is even.\\
Let $j_1 \in \Omega_1$ with $n_{2j_0+1}^2<m_{2j_1}$. Then from
Proposition $\ref{L.20}$, there exists a $(15,2j_1,\frac{1}{2})$
exact pair $(w_1,f_1)$ in $\mathfrak{X}_0$, with $w_1 \in <y_n:n
\in \N>$.
\\ Let $j_2=\sigma(f_1)$ and $F_2 \subset L$ with
$\#F_2=n_{2j_2}$ and $F_2>maxsupp(w_1)$.\\
We set $w_2=\frac{m_{2j_2}}{n_{2j_2}}\sum\limits_{i \in F_2}
e_i$ and
$f_2=\frac{1}{m_{2j_2}}\sum\limits_{i \in F_2} e_i^*$.
Then the pair $(w_2,f_2)$ is a $(15,2j_2,\frac{1}{2})$
exact pair in $\mathfrak{X}_0$.
\\ In this way we inductively construct the
$(15,2j_0+1,\frac{1}{2})$ dependent sequence
$(w_k,f_k)_{i=1}^{n_{2j_0+1}}$ in $\mathfrak{X}_0$.\\
From Proposition $\ref{L.2}$ we get that
\[\|\frac{1}{n_{2j_0+1}} \sum\limits_{k=1}^{n_{2j_0+1}}
(-1)^{k}w_k\|_{W_0} \le \frac{120}{m_{2j_0+1}^2}.\]
Hence
\begin{equation}\label{L.39}
\|T(\frac{1}{n_{2j_0+1}} \sum\limits_{k=1}^{n_{2j_0+1}}
(-1)^{k}w_k)\|_{W_0} \le \frac{120\|T\|}{m_{2j_0+1}^2}.
\end{equation}
We set \[A=\{k\in \{1,\ldots,n_{2j_0+1}\}:k \mbox{~even}\}
\mbox{~and~}
B=\{k\in \{1,\ldots,n_{2j_0+1}\}:k \mbox{~odd} \}.\]
Then
\begin{equation}\label{L.40}
\|\frac{1}{n_{2j_0+1}}\cdot\sum\limits_{k \in B}w_k\|_{W_0}
\ge
(\frac{1}{m_{2j_0+1}}\sum\limits_{k=1}^{n_{2j_0+1}} f_k)
(\frac{1}{n_{2j_0+1}} \cdot\sum\limits_{k \in B}w_k) \ge
\frac{1}{4m_{2j_0+1}}.
\end{equation}
and inequality $(\ref{L.33})$ yields that
\begin{equation}\label{L.29}
\|T(w_k)\|_{W_0}\le \frac{m_{2j_k}}{n_{2j_k}}d_2\le
\frac{1}{m_{2j_0+1}^2}d_2\le \frac{m}{8m_{2j_0+1}}
\mbox{~for every~} k \in A.
\end{equation}
Therefore from $(\ref{L.32}), (\ref{L.39}), (\ref{L.40}),
(\ref{L.29})$ and the triangle inequality we get that
\[\frac{120 \|T\|}{m_{2j_0+1}^2}\ge
m\|\frac{1}{n_{2j_0+1}}
\cdot\sum\limits_{k\in B} w_k\|_{W_0}-
\|T(\frac{1}{n_{2j_0+1}} \cdot
\sum\limits_{k\in A} w_k)\|_{W_0}\ge \frac{m}{4m_{2j_0+1}}-
\frac{m}{8m_{2j_0+1}}=\frac{m}{8m_{2j_0+1}}\]
which contradicts to the choice of $j_0 \in \N$.\\

\begin{Theorem}\label{L.22}
Let $T:\mathfrak{X}_0 \longrightarrow \mathfrak{X}_0$  be a
bounded linear operator. Then $T=\lambda I+S$, where $\lambda \in
\R,I$ the identity operator on $\mathfrak{X}_0$ and $S$ be a
bounded linear operator on $\mathfrak{X}_0$ such that whenever $Y$
is a subspace of $\mathfrak{X}_0$ with $S_{|Y}:Y \longrightarrow
S(Y)$ is an isomorphism it follows that the subspace $Y$ is
embedded isomorphically into $c_0(\N)$.
\end{Theorem}
\noindent{\bf{Proof:}} From Lemma $\ref{L.27}$ it follows that
there exist an $L\in [\N]$ and a sequence $(\lambda_n)_{n \in L}$
of real numbers such that $\lim\limits_{n \in L}
\|T(e_n)-\lambda_n e_n\|_{G_0}=0$. We observe that the sequence
$(\lambda_n)_{n \in L}$ is bounded and consequently there exists
$M \in [L]$ and $\lambda \in \R$ such that $\lim\limits_{n \in M}
\lambda_n=\lambda$. It is not hard to see that $\lim\limits_{n \in
M} \|T(e_n)-\lambda e_n\|_{G_0}=0$. Setting $S=T-\lambda I$, if
$S$ restricted in some subspace is an isomorphism, then from Lemma
$\ref{L.28}$ and Proposition $\ref{L.11}$ it follows the
conclusion.

\begin{Remark}
In a forthcoming paper we will present some variants of the space
$\mathfrak{X}_0$. More precisely for $1\le p<\infty$ we construct
an $\ell^p$ saturated Banach space $\mathfrak{X}_p$ such that for
every $p>1$ the space $\mathfrak{X}_p$ is reflexive and satisfies
tightness conditions similar to the corresponding ones of
$\mathfrak{X}_0$.
\end{Remark}

\end{document}